\documentclass[a4,11pt]{article}
\usepackage{color}
\usepackage{graphicx,graphics,epsfig}
\usepackage{amsmath,amsthm,amssymb}
\usepackage{array}

\def\XXint#1#2#3{{\setbox0=\hbox{$#1{#2#3}{\int}$ }
\vcenter{\hbox{$#2#3$ }}\kern-.5\wd0}}

\allowdisplaybreaks[2]

\newcommand{\R}{\mathbb{R}}
\newcommand{\N}{\mathbb{N}}

\begin{document}
\date{}
\title{A multi-moment scheme for the two dimensional Maxwell's equations.}

\author{Kazufumi Ito
\thanks{Center for Research in Scientific Computation \& Department of Mathematics,
North Carolina State University,
Raleigh,
NC 27695, USA.
({\tt kito@math.ncsu.edu})}
\and
Tomoya Takeuchi
\thanks{Center for Research in Scientific Computation,
North Carolina State University,
Raleigh,
NC 27695, USA.
({\tt tntakeuc@ncsu.edu})}
}

\maketitle
\begin{abstract}
We develop a numerical scheme for solving time-domain Maxwell's
equation.
The method is motivated by CIP method which uses function values and its derivatives as unknown variables.
The proposed scheme is developed by using the Poisson formula for the wave equation.
It is fully explicit space and
time integration method with higher order accuracy and CFL number
being one. The bi-cubic interpolation is used for the solution
profile to attain the resolution.
It preserves sharp profiles very accurately
without any smearing and distortion due to
the exact time integration and high resolution approximation. The
stability and numerical accuracy are investigated.
\end{abstract}

\section{Introduction}

Multi-moment methods for time dependent differential equations aim to increase the
accuracy of the numerical solution, and to lower the dispersive and
dissipative errors in the numerical solution. The most
distinguishing characteristic of the method is that more than one
moment per grid point or per cell, such as function value, its
derivatives and the integral over the cell etc, are considered as
unknown variables, and they are simultaneously updated by coupling the differential equation and derived differential equations for the derivatives.
A Hermite polynomial is usually defined on each cell using such quantities to interpolate the
numerical solution, and is used to update the solutions.
Such polynomial interpolation, being defined on each cell
individually, reduces the numerical stencil width. The compactness
of the stencil makes it feasible to handle the boundary condition
and the interface condition numerically.

Such schemes were first proposed for one dimensional hyperbolic
equations by Van-Leer \cite{Leer-Towaulticonsdiff:77} in a framework
of a finite difference method, and by Takewaki, Nishiguchi and Yabe
\cite{Takewaki+NishiguchiETAL-Cubiintepseumeth:85}
on the basis of the characteristic equation for the first order PDE.
The latter method is referred to as CIP method.
Their works triggered off development of multi moment methods, and
it is increasingly becoming an active area of research and have even
been applied to various equations.
There is a vast literature  developing CIP methods and the following is a partial list
of papers: nonlinear hyperbolic
equations \cite{Yabe+Aoki-univsolvhypeequa:91}, multi dimensional
hyperbolic equations
\cite{Yabe-univcubiintesolv:91,Yabe+IshikawaETAL-univsolvhypeequa:91},
a multi-dimensional the Maxwell's equations
\cite{Ogata+YabeETAL-AccuNumeScheMaxw:06}, a numerical simulation
for solid, liquid and gas
\cite{OGATA+KAWAGUCHI-NumeInveRareFlow:11,Yabe+TakizawaETAL-Chalunivsolvsoli:05,
Yabe-Univsimusolusoli:02}, a new mesh system applicable to
non-orthogonal coordinate system
\cite{Yabe+MizoeETAL-Highschewithmeth:04}, a variant of CIP method
\cite{Aoki-Intediffoper:97}.

CIP method uses the
cubic interpolation constructed via solution values and its
derivative at two end points of a cell to approximate the solution
in the cell.
Explicit time integration formulas for both the exact
solution and its derivative of a one dimensional transport equation with the constant
velocity field are coupled with the cubic interpolation to obtain an
explicit time integration numerical scheme (CIP scheme).
In \cite{Ito+Takeuchi-METHHYPESYSTWITH:11} we have developed and analyzed a CIP scheme for
one dimensional hyperbolic equations with variable and discontinuous coefficients.

Most of the proposed CIP methods so far, when it is applied to higher dimensional equations, relies on the Strang splitting technique: It reduces the equation of higher dimension to a sequence
of simpler one dimensional equations,  repeatedly applies a CIP method
for the reduced one dimensional equations.

In this paper we develop a multi-moment method for two dimensional Maxwell's equations, which does not employ the directional splitting technique at all.
We apply the exact integration method in time by the Poisson formula and use the
bi-cubic interpolation. We refer \cite{Cohen+Joly-Fourordeschehete:90} \cite{LeVeque-Finivolumethhype:02} \cite{Hesthaven+Warburton-NodaHighMethUnst:02},\cite{Quarteroni-Domadecomethsyst:90},
\cite{Cockburn+KarniadakisETAL:00} for numerical integration methods for the Maxwell and acoustic wave equations in general.
The original second order Yee scheme using the time and space staggered grid is
developed in \cite{Yee:66} and a fourth order time and space variant of
Yee scheme is developed in \cite{Deveze+BeaulieuETAL:92,Gustafsson+Mossberg-Timecomphighorde:04}.

Our contributions are as follows. The bi-cubic interpolation is
combined with the Poisson formula to develop a fully explicit time
and space numerical method for Maxwell, acoustic wave equations as
well as the second order wave equation. We analyze the von-Neumann
stability of the proposed method and we establish CFL number is one
for the proposed method. The one step of the method involves
updating four moments at each grid point and the symbolic
formulation and the symmetry of the update is used to reduce
operation counts. 
The method offers highly
efficient to attain the desired accuracy due to the relaxed CFL
number limitation and accurate space resolution by the bi-cubic
profile. The method is compared with the fourth order time-space Yee's scheme
in terms of numerical accuracy at nodes and required operation counts given accuracy.
The numerical convergence rate test shows the method is nearly fourth order and
the operation counts are comparable with those for the fourth order Yee's scheme
in conventional computing. With distributed computing implementation
the method becomes much efficient. Also, the method has the build-in bi-cubic interpolation
and thus provides sub-grid resolutions at each cell.

If we use the method with CFL$=1$, as shown in Section
\ref{sec:numericaltest} the method preserves sharp profiles in the solution
very accurately without any smearing and distortion. The stability
and numerical accuracy are analyzed. Also, the method computes
directly the physical quantities e.g., current and electric field
gradient, very accurately. The building block of our method is the
exact integration rule for the Poisson formula against the
polynomials. Thus, one can use various polynomial approximations
locally, including the one using the solution values only.
In the paper we only implement the methods for the periodic boundary
condition but it can be extended to the various boundary
conditions and class of absorbing boundary conditions.

Also, we can also extend the interface treatment in \cite{Li+Ito-immeintemeth:06} for piecewise cubic interpolation at the material discontinuity and develop the immersed interface method for discontinuous media.

An outline of our presentation is: in Section 2 a CIP scheme is proposed, in Section 3 the stability and error
analysis is presented. Finally in Section 4 we present our numerical tests and numerical convergent rate.

\section{Derivation of the multi-moment scheme}\label{sec:multi}
Consider the two dimensional (TE) Maxwell's equation for magnetic
field $H=(0,0,H_z)$ and electric field $E=(E_x,E_y,0)$:
\begin{equation} \label{Max}
\frac{\partial E_x}{\partial t} = \frac{1}{\varepsilon}\frac{\partial H_z}{\partial y},\quad
\frac{\partial E_y}{\partial t}= -\frac{1}{\varepsilon}\frac{\partial H_z}{\partial x},\quad
\frac{\partial H_z}{\partial t} = \frac{1}{\mu}\left(\frac{\partial E_x}{\partial y}
-\frac{\partial E_y}{\partial x}\right)
\end{equation}
where the material coefficients $\epsilon,\;\mu$ are constants. Let
$c=\frac{1}{\sqrt{\varepsilon\mu}}$ be the speed of light. Our
method uses the fact that \eqref{Max} is equivalent to second order
wave equations:
\begin{equation}\label{equ:maxwell}
  \partial^2_tH-c^2\Delta H_z = 0, \quad
  \partial^2_tE_x-c^2\Delta E_x = 0, \quad
  \partial^2_tE_y-c^2\Delta E_y = 0.
\end{equation}
under the assumption that div$(E_x,E_y)=0$.
Similarly, the two dimensional acoustic wave equation for
pressure $p$ and velocity $v=(v_x,v_y)$:
\begin{equation*} 
\rho \frac{\partial v_x}{\partial t} = \frac{\partial p}{\partial x},\quad
\rho \frac{\partial v_y}{\partial t} = \frac{\partial p}{\partial y},\quad
\frac{\partial p}{\partial t} = K\left(\frac{\partial v_x}{\partial
x}+\frac{\partial v_y}{\partial y}\right),
\end{equation*}
can be treated in exactly the same manner.

\subsection{Poisson's formula and the multi-moment scheme}

Let $u(t,x_1,x_2)$ be a solution of the wave equation:
\begin{equation*} \label{wave}
\partial^2_t u-c^2\Delta u = 0,\mbox{  with  } u(0)=g,\;\; u_t(0)=h.
\end{equation*}
By the Poisson's formula \cite{Evans-Partdiffequa:98}, we have
$$
u(t,x) =\frac{1}{2\pi ct}\int_{B(x,ct)}\frac{g(\xi) +
\nabla g(\xi)\cdot (\xi-x) + t
h(\xi)}{((ct)^2-|\xi-x|^2)^{\frac{1}{2}}}\,d\xi
$$
Using change of variable $\xi -x$ to $\xi$, the solution $u(t+\Delta t,x)$ is
{\small
\begin{equation*} \label{update}
\begin{array}{ll}
u(t+\Delta t, x)&=\displaystyle \frac{1}{2\pi c\Delta
t}\int_{B}\frac{u(t,x+\xi) + \xi \cdot
\nabla u(t,x+\xi)}{((c\Delta t)^2-|\xi|^2)^{\frac{1}{2}}}d\xi
d\eta +\frac{1}{2\pi c}\int_{B}\frac{\partial_t u(t,x+\xi)} {((c\Delta
t)^2-|\xi|^2)^{\frac{1}{2}}}d\xi,
\\ \\
& =L((1 + \xi\cdot\nabla) u(t,x+\cdot)|B ) + \Delta t
L(\partial_t u(t,x+\cdot))|B),
\end{array} \end{equation*}
}
where $B=B(0,c\Delta)$ and
\begin{equation*} \label{Lu}
L[u(t,\cdot) | B] =\frac{1}{2\pi c\Delta t}\int_{B}\frac{u(t,\xi)}{((c\Delta t)^2-|\xi|^2)^{\frac{1}{2}}}d\xi,
\end{equation*}
for a function $u(t,\xi)$ and $\xi=(\xi_1,\xi_2)$. \\

The derivatives $\partial_{x}^{\alpha}u(t,x)$, $\alpha=(\alpha_1,\alpha_2)$, $\alpha_1,\alpha_2 \in \mathbb{N}$, also satisfy the wave equation
\begin{align*}
\partial^2_t (\partial_{x}^{\alpha}u)-c^2\Delta (\partial_{x}^{\alpha}u) = 0,
\end{align*}
as long as the solution is smooth, and thus the higher order derivatives of the solution is advanced via
\begin{equation*}
\partial_{x}^{\alpha}u(t+\Delta t,x)
=L[(1+\xi\cdot \nabla) \partial_{\xi}^{\alpha}u(t,x+\cdot)|B]
+\Delta t L[\partial_{\xi}^{\alpha}\partial_tu(t,x+\cdot)|B].
\end{equation*}
We obtain the exact integration formula for the solutions and
its derivatives of $H_z(t,x)$, $E_x(t,x)$ and $E_y(t,x)$ at $x=0$:
For $\alpha =(0,0)$ and $\alpha \in \mathbb{N}^2$,
\begin{equation} \label{Mult}
\begin{array}{ll}
&\partial_{x}^{\alpha}H_z(t+\Delta t,x) =
  L[ (1 + \xi\cdot \nabla ) \partial_{\xi}^{\alpha}H_z(t,x+\cdot)|B]+
\Delta t L[\partial_{\xi}^{\alpha}\partial_tH_z(t,x+\cdot)|B]  \\ \\
 &=L[ (1 + \xi\cdot \nabla ) \partial_{\xi}^{\alpha}H_z(t,x+\cdot)|B]
 + \frac{\Delta t}{\mu}
 L[\partial_{\xi}^{\alpha} (\partial_{\xi_2} E_x(t,x+\cdot)-\partial_{\xi_1} E_y(t,x+\cdot) )|B],
\\ \\
&\partial_{x}^{\alpha}E_x(t+\Delta t,x) = L[ (1 + \xi\cdot \nabla )
\partial_{\xi}^{\alpha}E_{x}(t,x+\cdot)|B] +
\frac{\Delta
t}{\varepsilon}L[\partial_{\xi_2}\partial_{\xi}^{\alpha}H_z(t,x+\cdot)|B],
\\ \\
&\partial_{x}^{\alpha}E_y(t+\Delta t,x)
=L[ (1 + \xi\cdot \nabla ) \partial_{\xi}^{\alpha}E_{y}(x+\cdot)|B]-\frac{\Delta
t}{\varepsilon}
L[\partial_{\xi_1} \partial_{\xi}^{\alpha}H_z(t,x+\cdot)|B].
\end{array} \end{equation}
Here we use \eqref{equ:maxwell} to exchange the time derivative and spatial derivatives.

Let us define a grid of points in the $(t,x)$ space. Let $\Delta t$ and $\Delta x$ be positive numbers. The grid is the set of points $(t_n,x_{ij})=(n\Delta t, i\Delta, j\Delta)$ for arbitrary integers $(n,i,j)$. 
We let $\partial_{x}^{\alpha}{H_z}^{n}_{ij}$, $\partial_{x}^{\alpha}{E_x}^{n}_{ij}$ and $\partial_{x}^{\alpha}{E_y}^{n}_{ij}$ stand  respectively for the approximation to the solution
$\partial_{x}^{\alpha}{H_z}(t_n,x_i,y_j)$, $\partial_{x}^{\alpha}E_x(t_n,x_i,y_j)$ and  $\partial_{x}^{\alpha}E_y(t_n,x_i,y_j)$ for $\alpha \in(0,0)\cup\mathbb{N}$.

The basis idea of the multi-moment scheme is to define a higher order polynomials $P(x,y)$ on each cell $[x_{i},x_{i+1}]\times[y_j,y_{j+1}]$ using grid values including spatial derivatives at four corners of the cell, and substitute them to the exact time integration formula \eqref{Mult}: We evaluate the integrals of the polynomials over the ball $B(0,c\Delta t)$. Henceforth we assume that $c\Delta t \le \Delta x $, and thus the four polynomials are involved in the integration over the ball, and thus the method uses variables at 9 nearest grid points.

We can derive various multi-moment schemes on the basis of \eqref{Mult}. The resulted scheme depends on the number of unknown variables we employ at each grid and the order of interpolation polynomials; for instance, if we take the grid values, the firs-order derivatives and their second order mixed derivatives as unknown variables, we use the bi-cubic Hermite interpolation, known as Boger-Fox-Schmit element in finite element methods,
\begin{equation*}
\sum_{k=0}^3\sum_{\ell=0}^3 c_{k,\ell} x^ky^{\ell}.
\end{equation*}
The coefficients of the polynomial are defined by the usual interpolation condition. The resulted scheme is written in terms of the bi-cubic polynomial. Let ${H_z}_{,ij}$, ${E_x}_{,ij}$ and ${E_y}_{,ij}$ denote the bi-cubic polynomials defined in the cell $[x_i,x_{i+1}]\times [y_j,y_{j+1}]$ by the interpolation condition:
\begin{equation*} \label{bi-cubic}
\partial_x^\alpha{H_z}_{,ij}(x_{\ell m})=\partial_x^\alpha {H_z}^n_{,\ell m}, \quad
\partial_x^\alpha{E_x}_{,ij}(x_{\ell m})=\partial_x^\alpha {E_x}^n_{,\ell m}, \quad
\partial_x^\alpha{E_y}_{,ij}(x_{\ell m})=\partial_x^\alpha {E_y}^n_{,\ell m},
\end{equation*}
for $\alpha \in\{(0,0),(1,0),(0,1),(1,1)\}$ and $x_{\ell m}\in\{x_{ij},x_{i-1,j},x_{i-1,j-1},x_{i,j-1}\}$. Let us number four cells surrounding a gird $x_{ij}$ counter clockwise; $C_1=[x_i,x_{i+1}]\times [y_j,y_{j+1}]$, $C_2=[x_{i-1},x_{i}]\times [y_j,y_{j+1}]$, $C_3=[x_{i-1},x_{i}]\times [y_{j-1},y_{j}]$ and $C_4=[x_i,x_{i+1}]\times [y_{j-1},y_{j}]$. We also number the polynomial defined on each cell accordingly, i.e., ${H_z}_{,ij} = {H_z}_{,1}$, etc.
We let $B_k$ stand for
\begin{equation*}
\begin{array}{ll}
B_1=B(0,c\Delta t)\cap\{x\ge 0\}\cap \{y\ge 0\},
& B_2=B(0,c\Delta t)\cap\{x\le 0\}\cap \{y\ge 0\},\\
B_3=B(0,c\Delta t)\cap\{x\le 0\}\cap \{y\le 0\},
& B_4=B(0,c\Delta t)\cap\{x\ge 0\}\cap \{y\le 0\}.
\end{array}
\end{equation*}
With these notation, we obtain the CIP scheme:
{\small
 \begin{equation} \label{MultPolyForm}
\left\{\begin{array}{l}
\partial_{x}^{\alpha}{H_z}^{n+1}_{ij} =
\sum_{k=1}^4 L[ (1 + \xi\cdot \nabla ) \partial_{\xi}^{\alpha}{H_z}_{,k}(x_{ij}+\cdot) |B_k]\\
\\\quad  + \frac{\Delta t}{\mu}
\sum_{k=1}^4 L[\partial_{\xi}^{\alpha} (\partial_{\xi_2} {E_x}_{,k}(x_{ij}+\cdot)-\partial_{\xi_1} {E_y}_{,k}(x_{ij}+\cdot) )|B_k],
\\ \\
\partial_{x}^{\alpha}{E_x}^{n+1}_{ij} = \sum_{k=1}^4L[ (1 + \xi\cdot \nabla )
\partial_{\xi}^{\alpha}{E_{x}}_{,k}(x_{ij}+\cdot)|B_k] +
\frac{\Delta
t}{\varepsilon}\sum_{k=1}^4L[\partial_{\xi_2}\partial_{\xi}^{\alpha}{H_z}_{,k}(x_{ij}+\cdot)|B_k],
\\ \\
\partial_{x}^{\alpha}{E_y}^{n+1}_{ij}
=\sum_{k=1}^4L[ (1 + \xi\cdot \nabla ) \partial_{\xi}^{\alpha}{E_{y}}_{,k}(x_{ij}+\cdot)|B_k]-\frac{\Delta
t}{\varepsilon}
\sum_{k=1}^4L[\partial_{\xi_1} \partial_{\xi}^{\alpha}{H_z}_{,k}(x_{ij}+\cdot)|B_k],
\end{array} \right.\end{equation}
}
 for $\alpha=(0,0)$ (function value update), $\alpha=(1,0)$ ($x_1$ derivative update), $\alpha=(0,1)$ ($x_2$ derivative update) and $\alpha=(1,1)$ ($x_1,x_2$ second order mixed derivative update). As detailed in the following sections, \eqref{MultPolyForm} develops the moments update at time step $t_{n+1}$ based on the 9 nearest grid moments at time step $t_n$, i.e., update \eqref{alg0}.

One can reduce the number of unknowns at a grid point; for example, one uses the grid values, the firs-order derivatives as unknowns. The number of unknowns to be determined at each grid point becomes 9: each component $H_z,E_x,E_y$ has 3 unknowns at a grid point. Possible choices for the interpolation are
\begin{align*}
 \sum_{k=0}^3\sum_{\ell=0}^k c_{k,\ell} x^{k-\ell}y^k + c_{3,4}x^2y^3 + c_{4,3}x^3y^2,\quad \mbox{or}\quad
\sum_{k=0}^3\sum_{\ell=0}^k c_{k,\ell} x^{k-\ell}y^k + c_{4,1}x^3y + c_{1,4}xy^3.
\end{align*}
The coefficients $c_{k,\ell}$ are determined by using the grid values, the firs-order  derivatives at the four corners of a cell. The second mixed derivatives being not used, the resulted schemes have less complexity than the bi-cubic Hermite polynomial based scheme, however, they produce less accurate numerical solution, and the CFL number is less than 1. As for the other choice, we consider the bi-liner interpolation:
\begin{align*}
c_{0,0} + c_{1,0}x + c_{0,1}y + c_{1,1}xy.
\end{align*}
We then obtain a derivative free nine point scheme.
\subsection{Bi-cubic interpolation and the integration}

In this section, we examine the details for computing the integrals in \eqref{MultPolyForm}
when the bi-cubic interpolation is used for the interpolation, i.e., the number of unknowns for  $H_z$, $E_x$ and $E_y$ is 4 respectively at a grid point; the grid value, the firs-order derivatives and the second order mixed derivative.
\subsubsection{Notation}
Let us introduce some notation.
For the numerical quantities $f_{i,j}$, $\partial_{x_1}f_{i,j}$,
$\partial_{x_2}f_{i,j}$ and $\partial^2_{x_1x_2}f_{i,j}$ given
at the node $x_{ij}$,
we define a $4 \times 1$ vector (a \textit{multi-moment vector}):
\begin{equation*}
\mathbf{f}_{i,j} = \left[ f_{i,j}, \partial_{x_1}f_{i,j},
\partial_{x_2}f_{i,j},\partial^2_{x_1,x_2}f_{i,j}\right]^\top
 \in \mathbb{R}^4,
\end{equation*}
Let $f^{i,j}$ denote a vector composed of the multi-moment vector assigned at the four corner of the cell $C_{i,j}=[x_i,x_{i+1}] \times [y_j,y_{j+1}]$:
\begin{equation*}
f^{i,j} = \left[
       \mathbf{f}_{i,j}^\top ,
       \mathbf{f}_{i+1,j}^\top  ,
       \mathbf{f}_{i+1,j+1}^\top  ,
       \mathbf{f}_{i,j+1}^\top ,
  \right]^\top \in \mathbb{R}^{16}.
\end{equation*}

Denote $\frac{x-x_{ij}}{d_{ij}}:= \left(\frac{x_1-{x_1}_{,i}}{d_{1,i}},  \frac{x_2-{x_2}_{,i}}{ d_{2,j}}\right)$ where $d_{1,i}=x_{1,i+1}-x_{1,i} $ and $d_{2,j}= x_{2,j+1}-x_{2,j}$.  We construct a bi-cubic polynomial
\begin{equation*}
F_{i,j}(x)= \sum_{k=0}^3\sum_{\ell=0}^3 q_{k,\ell}
\left(\frac{x_1-{x_1}_{,i}}{d_{1,i}}\right)^{k}
\left(\frac{x_2-{x_2}_{,j}}{d_{2,j}}\right)^\ell
= e\left(\frac{x-x_{ij}}{d_{ij}}\right) q^{i,j},  \label{bcub0}
\end{equation*}
on the cell $C_{i,j}$, where the coefficient vector $q^{i,j}$ are ordered as
\begin{equation}\label{q}
 q^{i,j}:=(q_{0,0},q_{0,1},q_{0,2},q_{0,3},q_{1,0},q_{1,1},q_{1,2},q_{1,3},q_{2,0},q_{2,1},
 q_{2,2},q_{2,3},q_{3,0},q_{3,1},q_{3,2},q_{3,3})^\top.
\end{equation}
and $e(x)$ denotes $1\times 16$ row vectors
\begin{equation*}
  e(x) = \left[1,x_1 ,x_1^2,x_1^3
  \right]
  \otimes \left[1,x_2 ,x_2^2,x_2^3
  \right],
\end{equation*}
i.e., the components of $e(x)$ are
\begin{equation*}\label{e}
e_{\ell,m}(x_1,x_2)=x_1^\ell x_2^m,
\end{equation*}
and are ordered as in \eqref{q}. \\
The coefficient $q^{i,j}$ of the bi-cubic polynomial $F_{i,j}$ is determined by 16 interpolation
conditions at four corners $x_{i,j}$, $x_{i+1,j}$, $x_{i,j+1}$ and $x_{i+1,j+1}$ of
the cell:
\begin{equation*} \label{bcub1}
\begin{smallmatrix}
F_{i,j}(x_{i,j} )=f_{i,j}, & \partial_{x_1}F_{i,j}(x_{i,j}) =
\partial_{x_1}f_{i,j}, &
\partial_{x_2}F_{i,j}(x_{i,j})=\partial_{x_2}f_{i,j},
&\partial^2_{x_1x_2}F_{i,j}(x_{i,j})=\partial^2_{x_1x_2}f_{i,j}, \\ \\
F_{i,j}(x_{i+1,j} )=f_{i+1,j}, & \partial_{x_1}F_{i,j}(x_{i+1,j}) =
\partial_{x_1}f_{i+1,j}, &
\partial_{x_2}F_{i,j}(x_{i+1,j})= \partial_{x_2}f_{i+1,j},
& \partial^2_{x_1x_2}F_{i,j}(x_{i+1,j}) =\partial^2_{x_1x_2}f_{i+1,j}, \\ \\
F_{i,j}(x_{i+1,j+1} )=f_{i+1,j+1}, & \partial_{x_1}F_{i,j}(x_{i+1,j+1}) =
\partial_{x_1}f_{i+1,j+1}, &
\partial_{x_2}F_{i,j}(x_{i+1,j+1})=\partial_{x_2}f_{i+1,j+1},
& \partial^2_{x_1x_2}F_{i,j}(x_{i+1,j+1}) = \partial^2_{x_1x_2}f_{i+1,j+1}, \\ \\
F_{i,j}(x_{i,j+1})=f_{i,j+1}, & \partial_{x_1}F_{i,j}(x_{i,j+1}) =
\partial_{x_1}f_{i,j+1}, &
\partial_{x_2}F_{i,j}(x_{i,j+1})=\partial_{x_2}f_{i,j+1}, &
\partial^2_{x_1x_2}F_{i,j}(x_{i,j+1}) =
\partial^2_{x_1x_2}f_{i,j+1}.
\end{smallmatrix}
\end{equation*}

We obtain then the coefficients $q^{i,j}$ of $F_{i,j}$:
\[
q^{i,j}=Q R_{ij}f^{i,j},
\]
where $Q = [Q_1, Q_2, Q_3, Q_4]$ is the interpolation matrix:
\begin{equation*}
Q_1 =\left[\begin{smallmatrix}
 1 & 0 & 0 & 0 \\
 0 & 0 & 1 & 0 \\
 -3 & 0 & -2 & 0 \\
 2 & 0 & 1 & 0 \\
 0 & 1 & 0 & 0 \\
 0 & 0 & 0 & 1 \\
 0 & -3 & 0 & -2 \\
 0 & 2 & 0 & 1 \\
 -3 & -2 & 0 & 0 \\
 0 & 0 & -3 & -2 \\
 9 & 6 & 6 & 4 \\
 -6 & -4 & -3 & -2 \\
 2 & 1 & 0 & 0 \\
 0 & 0 & 2 & 1 \\
 -6 & -3 & -4 & -2 \\
 4 & 2 & 2 & 1
  \end{smallmatrix}\right],
Q_2 = \left[\begin{smallmatrix}
 0 & 0 & 0 & 0 \\
 0 & 0 & 0 & 0 \\
 0 & 0 & 0 & 0 \\
 0 & 0 & 0 & 0 \\
 0 & 0 & 0 & 0 \\
 0 & 0 & 0 & 0 \\
 0 & 0 & 0 & 0 \\
 0 & 0 & 0 & 0 \\
 3 & -1 & 0 & 0 \\
 0 & 0 & 3 & -1 \\
 -9 & 3 & -6 & 2 \\
 6 & -2 & 3 & -1 \\
 -2 & 1 & 0 & 0 \\
 0 & 0 & -2 & 1 \\
 6 & -3 & 4 & -2 \\
 -4 & 2 & -2 & 1
  \end{smallmatrix}\right] ,
  Q_3 =\left[\begin{smallmatrix}
 0 & 0 & 0 & 0 \\
 0 & 0 & 0 & 0 \\
 0 & 0 & 0 & 0 \\
 0 & 0 & 0 & 0 \\
 0 & 0 & 0 & 0 \\
 0 & 0 & 0 & 0 \\
 0 & 0 & 0 & 0 \\
 0 & 0 & 0 & 0 \\
 0 & 0 & 0 & 0 \\
 0 & 0 & 0 & 0 \\
 9 & -3 & -3 & 1 \\
 -6 & 2 & 3 & -1 \\
 0 & 0 & 0 & 0 \\
 0 & 0 & 0 & 0 \\
 -6 & 3 & 2 & -1 \\
 4 & -2 & -2 & 1
  \end{smallmatrix}\right],
  Q_4 = \left[  \begin{smallmatrix}
 0 & 0 & 0 & 0 \\
 0 & 0 & 0 & 0 \\
 3 & 0 & -1 & 0 \\
 -2 & 0 & 1 & 0 \\
 0 & 0 & 0 & 0 \\
 0 & 0 & 0 & 0 \\
 0 & 3 & 0 & -1 \\
 0 & -2 & 0 & 1 \\
 0 & 0 & 0 & 0 \\
 0 & 0 & 0 & 0 \\
 -9 & -6 & 3 & 2 \\
 6 & 4 & -3 & -2 \\
 0 & 0 & 0 & 0 \\
 0 & 0 & 0 & 0 \\
 6 & 3 & -2 & -1 \\
 -4 & -2 & 2 & 1
  \end{smallmatrix}\right].
\end{equation*}
And $R_{ij}$ is a tensor product of the 4 by 4 identity matrix $I$ and the diagonal matrix with the diagonal entries $[1, \Delta x_{1,i}, \Delta x_{2,j}, \Delta x_{1,i}\Delta x_{2,j}]$:
\begin{equation*}
  R_{ij}=I\otimes \left[
                      \begin{array}{cccc}
                        1 &  &  &  \\
                         & d_{1,i} & & \\
                         & & d_{2,j} & \\
                         &  &  & d_{1,i} d_{2,j} \\
                      \end{array}
                    \right].
\end{equation*}
Thus we obtain the bi-cubic polynomial:
\begin{align*}\label{bicubic}
F_{i,j}(x) = e\left(\frac{x-x_{ij}}{d_{ij}}\right)QR_{ij}f^{i,j} 
&= e\left(\frac{x-x_{ij}}{d_{ij}}\right)Q_1R_{i,j}\mathbf{f}_{i,j} + e\left(\frac{x-x_{ij}}{d_{ij}}\right)Q_2R_{ij}\mathbf{f}_{i+1,j}\\
&+e\left(\frac{x-x_{ij}}{d_{ij}}\right)Q_3R_{ij}\mathbf{f}_{i+1,j+1} +e\left(\frac{x-x_{ij}}{d_{ij}}\right)Q_4R_{ij}\mathbf{f}_{i,j+1}.\nonumber
\end{align*}
Next, let us introduce some matrices for basic operations.
For a cubic polynomial $p(x) = e_0(x)a$ for $x\in \mathbb{R}$, where $e_0(x) =
(1,x,x^2,x^3)
$
 and $a=(a_0, a_1, a_2, a_3)^\top$, we have
\begin{equation*}
   \frac{d}{dx}p(x)= e_0(x)Da,\quad
  x\frac{d}{dx}p(x) =e_0(x)M_xDa ,\quad
   p(\alpha x) = e_0(x)D_\alpha a,\quad
  p(x-s) = e_0(x) T_s a,
\end{equation*}
where{\footnotesize
\begin{align*}\label{D_etc}
  D = \left[
        \begin{array}{cccc}
          0 & 1 & 0 & 0 \\
          0 & 0 & 2 & 0 \\
          0 & 0 & 0 & 3 \\
          0 & 0 & 0 & 0 \\
        \end{array}
      \right],\quad
  M=\left[
      \begin{array}{cccc}
        0 & 0 & 0 & 0 \\
        1 & 0 & 0 & 0 \\
        0 & 1 & 0 & 0 \\
        0 & 0 & 1 & 0 \\
      \end{array}
    \right],
   D_\alpha=\left[
      \begin{array}{cccc}
        1 & 0 & 0 & 0 \\
        0 & \alpha & 0 & 0 \\
        0 & 0 & \alpha^2 & 0 \\
        0 & 0 & 0 & \alpha^3 \\
      \end{array}
    \right],\quad
      T_x=\left[
   \begin{array}{cccc}
1 & -x & x^2 & -x^3 \\
0 & 1 & -2x & 3x^2 \\
 0 & 0 & 1 & -3x \\
0 & 0 & 0 & 1 \\
\end{array}
   \right].
\end{align*}
}
Below, we will use the commutative properties:
\begin{equation}\label{property}
 D D_\alpha = {\alpha}D_\alpha D,\quad
  M D_\alpha ={\alpha}^{-1} D_\alpha M,\quad
  T_s D_\alpha = D_\alpha T_s.
\end{equation}
\subsubsection{Computation of $L$}
Now we express the integral in \eqref{MultPolyForm} in terms of the grid values.
We compute the integrals
\begin{equation*}
L(A F_{k}(x_{i,j}+\cdot)|B_k),
\end{equation*}
for $k=1,2,3,4$,
where $A$ denotes one of the operators
\begin{equation*}\label{Ts2}
\begin{array}{l}
(1+\xi\cdot \nabla)\partial^{\alpha_1}_{\xi_1} \partial^{\alpha_2}_{\xi_2},\quad
 \partial_{\xi_1}\partial^{\alpha_1}_{\xi_1} \partial^{\alpha_2}_{\xi_2},\quad \partial_{\xi_2}\partial^{\alpha_1}_{\xi_1} \partial^{\alpha_2}_{\xi_2},
\end{array}
\end{equation*}
for $\alpha_1=0,1$, $\alpha_2=0,1$,
and we renumber the polynomial; $F_1=F_{i,j}$, $F_2=F_{i-1,j}$, $F_3=F_{i-1,j-1}$, and $F_4=F_{i,j-1}$.
Let us denote the corresponding matrix representation for $A$ by $T_A$, i.e.,
\begin{align*}
 &T_{(1+\xi\cdot \nabla)\partial^{\alpha_1}_{\xi_1} \partial^{\alpha_2}_{\xi_2}}= (I\otimes I+ MD\otimes I + I \otimes MD) (D^{\alpha_1}\otimes D^{\alpha_2}), \\
&T_{\partial_{\xi_1}\partial^{\alpha_1}_{\xi_1} \partial^{\alpha_2}_{\xi_2}}=  (D\otimes I)  (D^{\alpha_1}\otimes D^{\alpha_2}),\quad
T_{\partial_{\xi_2}\partial^{\alpha_1}_{\xi_1} \partial^{\alpha_2}_{\xi_2}}=  (I\otimes D)  (D^{\alpha_1}\otimes D^{\alpha_2}).\nonumber
\end{align*}
For the compact expression, we also number the multi-moment vector accordingly, i.e.,
$f^{1}=f^{i,j}$, $f^{2}=f^{i-1,j}$, $f^{3}=f^{i-1,j-1}$, and $f^{4}=f^{i,j-1}$.
Then from the representations
\begin{align*}
&F_{1}(x_{i,j}+\xi) = e\left(\frac{\xi+x_{i,j}-x_{i,j}}{d_{ij}}\right) QR_{i,j} f^{1}=
e(\xi)\left(D_{\frac{1}{d_{1,i}}}\otimes D_{\frac{1}{d_{2,j}}}\right) QR_{i,j} f^{1}, \\
&F_{2}(x_{i,j}+\xi) = e(\xi)\left(D_{\frac{1}{d_{1,i-1}}}\otimes D_{\frac{1}{d_{2,j}}}\right)(T_{-1} \otimes I)QR_{i-1,j}f^{2}, \\
&F_{3}(x_{i,j}+\xi) = e(\xi) \left(D_{\frac{1}{d_{1,i-1}}}\otimes D_{\frac{1}{d_{2,j-1}}}\right)( T_{-1}\otimes T_{-1})QR_{i-1,j-1}f^{3},\\
&F_{4}(x_{i,j}+\xi) = e(\xi)\left(D_{\frac{1}{d_{1,i}}}\otimes D_{\frac{1}{d_{2,j-1}}}\right)(I\otimes T_{-1})QR_{i,j-1}f^{4},
\end{align*}
and hence we have
{\small
\begin{align*}\label{basic}
&L(A F_{1}(x_{i,j}+\cdot) | B_1)=L(e| B_1)\,T_A \left(D_{\frac{1}{d_{1,i}}}\otimes D_{\frac{1}{d_{2,j}}}\right) QR_{i,j} f^{1}, \\
&L(A F_{2}(x_{i,j}+\cdot)| B_2) =L(e| B_2) \,T_A \left(D_{\frac{1}{d_{1,i-1}}}\otimes D_{\frac{1}{d_{2,j}}}\right)(T_{-1} \otimes I)QR_{i-1,j}f^{2}\\
&L(A F_{3}(x_{i,j}+\cdot) | B_3)=L(e| B_3) \,T_A \left(D_{\frac{1}{d_{1,i-1}}}\otimes D_{\frac{1}{d_{2,j-1}}}\right)( T_{-1}\otimes T_{-1})QR_{i-1,j-1}f^{3}\\
&L(A F_{4}(x_{i,j}+\cdot) | B_4)=L(e| B_4) \,T_A \left(D_{\frac{1}{d_{1,i}}}\otimes D_{\frac{1}{d_{2,j-1}}}\right)(I\otimes T_{-1})QR_{i,j-1}f^{4}.
\end{align*}
}
Thus the computation of the integrals is reduced to the computations of
\begin{equation*} 
L(e_{\ell,m}| B_i)=\displaystyle \frac{1}{2\pi c\Delta t}
\int_{B_i}\frac{e_{\ell,m}(\xi)}{((c\Delta t)^2-|\xi|^2)^{\frac{1}{2}}}\,d\xi .
\end{equation*}
Using change of variable $\xi=(\xi_1,\xi_2)=c\Delta t\,
r\,(\cos\theta,\sin\theta)$, the integrals are evaluated as the function of $d_c:=c\Delta t$:
{\small
\begin{align*}
{\mathbf{d_{1,c}}}&:=L(e|B_1)
=\left[\frac{1}{4},\frac{{d_c}}{8},\frac{{d_c}^2}{12},\frac{{d_c}^3}{16},\frac{{d_c}}{8},\frac{{d_c}^2}{6 \pi },\frac{{d_c}^3}{32},\frac{{d_c}^4}{15 \pi },\frac{{d_c}^2}{12},\frac{{d_c}^3}{32},\frac{{d_c}^4}{60},\frac{{d_c}^5}{96},\frac{{d_c}^3}{16},\frac{{d_c}^4}{15 \pi },\frac{{d_c}^5}{96},\frac{2 {d_c}^6}{105 \pi }\right],\\
{\mathbf{{d}_{2,c}}}&:=L(e|B_2)
=\left[\frac{1}{4},\frac{{d_c}}{8},\frac{{d_c}^2}{12},\frac{{d_c}^3}{16},-\frac{{d_c}}{8},-\frac{{d_c}^2}{6 \pi },-\frac{{d_c}^3}{32},-\frac{{d_c}^4}{15 \pi },\frac{{d_c}^2}{12},\frac{{d_c}^3}{32},\frac{{d_c}^4}{60},\frac{{d_c}^5}{96},-\frac{{d_c}^3}{16},-\frac{{d_c}^4}{15 \pi },-\frac{{d_c}^5}{96},-\frac{2 {d_c}^6}{105 \pi }\right],\\
{\mathbf{{d}_{3,c}}}&:=L(e|B_3)
=\left[\frac{1}{4},-\frac{{d_c}}{8},\frac{{d_c}^2}{12},-\frac{{d_c}^3}{16},-\frac{{d_c}}{8},\frac{{d_c}^2}{6 \pi },-\frac{{d_c}^3}{32},\frac{{d_c}^4}{15 \pi },\frac{{d_c}^2}{12},-\frac{{d_c}^3}{32},\frac{{d_c}^4}{60},-\frac{{d_c}^5}{96},-\frac{{d_c}^3}{16},\frac{{d_c}^4}{15 \pi },-\frac{{d_c}^5}{96},\frac{2 {d_c}^6}{105 \pi }\right],\\
{\mathbf{{d}_{4,c}}}&:=L(e|B_4)
=\left[\frac{1}{4},-\frac{{d_c}}{8},\frac{{d_c}^2}{12},-\frac{{d_c}^3}{16},\frac{{d_c}}{8},-\frac{{d_c}^2}{6 \pi },\frac{{d_c}^3}{32},-\frac{{d_c}^4}{15 \pi },\frac{{d_c}^2}{12},-\frac{{d_c}^3}{32},\frac{{d_c}^4}{60},-\frac{{d_c}^5}{96},\frac{{d_c}^3}{16},-\frac{{d_c}^4}{15 \pi },\frac{{d_c}^5}{96},-\frac{2 {d_c}^6}{105 \pi }\right].
\end{align*}
}
Using these vectors, the integrations are expressed in terms of the vectors and the matrices, i.e., we obtain
{\small
\begin{align*}
&L(A F_{1}(x_{i,j}+\cdot) | B_1)={\mathbf{d_{1,c}}}\,T_A \left(D_{\frac{1}{d_{1,i}}}\otimes D_{\frac{1}{d_{2,j}}}\right) QR_{i,j} f^{1}, \\
&L(A F_{2}(x_{i,j}+\cdot)| B_2) ={\mathbf{d_{2,c}}}\,T_A \left(D_{\frac{1}{d_{1,i-1}}}\otimes D_{\frac{1}{d_{2,j}}}\right)(T_{-1} \otimes I)QR_{i-1,j}f^{2}\\
&L(A F_{3}(x_{i,j}+\cdot) | B_3)={\mathbf{d_{3,c}}} \,T_A \left(D_{\frac{1}{d_{1,i-1}}}\otimes D_{\frac{1}{d_{2,j-1}}}\right)( T_{-1}\otimes T_{-1})QR_{i-1,j-1}f^{3}\\
&L(A F_{4}(x_{i,j}+\cdot) | B_4)={\mathbf{d_{4,c}}}\,T_A \left(D_{\frac{1}{d_{1,i}}}\otimes D_{\frac{1}{d_{2,j-1}}}\right)(I\otimes T_{-1})QR_{i,j-1}f^{4}.
\end{align*}
}
From \eqref{property}, we have for $d_1,d_2 > 0$
{\small
\begin{align*}
 & T_{(1+\xi\cdot \nabla)\partial^{\alpha_1}_{\xi_1} \partial^{\alpha_2}_{\xi_2}} (D_{d_1^{-1}}\otimes D_{d_2^{-1}}) =
   \frac{1}{d_1^{\alpha_1} d_2^{\alpha_2}}(D_{d_1^{-1}}\otimes D_{d_2^{-1}})T_{(1+\xi\cdot \nabla)\partial^{\alpha_1}_{\xi_1} \partial^{\alpha_2}_{\xi_2}}
  \\
 &  T_{\partial_{\xi_1}\partial^{\alpha_1}_{\xi_1} \partial^{\alpha_2}_{\xi_2}}(D_{d_1^{-1}}\otimes D_{d_2^{-1}}) =\frac{1}{d_1^{\alpha_1} d_2^{\alpha_2}}\frac{1}{d_1}(D_{d_1^{-1}}\otimes D_{d_2^{-1}}) T_{\partial_{\xi_1}\partial^{\alpha_1}_{\xi_1} \partial^{\alpha_2}_{\xi_2}}  ,\\
 &  T_{\partial_{\xi_2}\partial^{\alpha_1}_{\xi_1} \partial^{\alpha_2}_{\xi_2}}(D_{d_1^{-1}}\otimes D_{d_2^{-1}}) =\frac{1}{d_1^{\alpha_1} d_2^{\alpha_2}}\frac{1}{d_2}(D_{d_1^{-1}}\otimes D_{d_2^{-1}}) T_{\partial_{\xi_2}\partial^{\alpha_1}_{\xi_1} \partial^{\alpha_2}_{\xi_2}}.
\end{align*}
}
We compute
{\small
\begin{align*}
&{\mathbf{d_{1,c}}}\left(D_{\frac{1}{d_{1,i}}}\otimes D_{\frac{1}{d_{2,j}}}\right)=\left[\frac{1}{4},\frac{\mu }{8},\frac{\mu ^2}{12},\frac{\mu ^3}{16},\frac{\lambda }{8},\frac{\lambda  \mu }{6 \pi },\frac{\lambda  \mu ^2}{32},\frac{\lambda  \mu ^3}{15 \pi },\frac{\lambda ^2}{12},\frac{\lambda ^2 \mu }{32},\frac{\lambda ^2 \mu ^2}{60},\frac{\lambda ^2 \mu ^3}{96},\frac{\lambda ^3}{16},\frac{\lambda ^3 \mu }{15 \pi },\frac{\lambda ^3 \mu ^2}{96},\frac{2 \lambda ^3 \mu ^3}{105 \pi }\right],
\end{align*}
}
where $\lambda = \frac{c\Delta t}{d_{1,i}}$ and $\mu = \frac{c\Delta t}{d_{2,j}}$.
We denote the right hand side by $ \Lambda(\lambda,\mu)$, i.e.,
{\small
\begin{equation*}
\Lambda(\lambda,\mu)=\left[\frac{1}{4},\frac{\mu }{8},\frac{\mu ^2}{12},\frac{\mu ^3}{16},\frac{\lambda }{8},\frac{\lambda  \mu }{6 \pi },\frac{\lambda  \mu ^2}{32},\frac{\lambda  \mu ^3}{15 \pi },\frac{\lambda ^2}{12},\frac{\lambda ^2 \mu }{32},\frac{\lambda ^2 \mu ^2}{60},\frac{\lambda ^2 \mu ^3}{96},\frac{\lambda ^3}{16},\frac{\lambda ^3 \mu }{15 \pi },\frac{\lambda ^3 \mu ^2}{96},\frac{2 \lambda ^3 \mu ^3}{105 \pi }\right].
\end{equation*}
}
Let $\lambda_1=\lambda_4=\frac{c\Delta t}{d_{1,i}}$, $\lambda_2=\lambda_3=\frac{c\Delta t}{d_{1,i-1}}$,
$\mu_1=\mu_2=\frac{c\Delta t}{d_{2,j}}$ and $\mu_3=\mu_4=\frac{c\Delta t}{d_{2,j-1}}$. Then
{\small
\begin{align*}
&{\mathbf{d_{1,c}}}\left(D_{\frac{1}{d_{1,i}}}\otimes D_{\frac{1}{d_{2,j}}}\right)=\Lambda(\lambda_1,\mu_1),\quad
{\mathbf{d_{2,c}}}\left(D_{\frac{1}{d_{1,i-1}}}\otimes D_{\frac{1}{d_{2,j}}}\right)=\Lambda(-\lambda_2,\mu_2),\\
&{\mathbf{d_{3,c}}}\left(D_{\frac{1}{d_{1,i-1}}}\otimes D_{\frac{1}{d_{2,j-1}}}\right)=\Lambda(-\lambda_3,-\mu_3),\quad
{\mathbf{d_{4,c}}}\left(D_{\frac{1}{d_{1,i}}}\otimes D_{\frac{1}{d_{2,j-1}}}\right)=\Lambda(\lambda_4,-\mu_4)
\end{align*}
}
Therefore we have for $A=(1+\xi\cdot \nabla)\partial^{\alpha_1}_{\xi_1} \partial^{\alpha_2}_{\xi_2}$,
{\small
\begin{align}\label{formula1}
&\sum_{k=1}^4 L(A F_{k}(x_{i,j}+\cdot) | B_k)=
\frac{\Lambda(\lambda_1,\mu_1)}{d_{1,i}^{\alpha_1} d_{2,j}^{\alpha_2}}\,T_{A}QR_{i,j} f^{1}+
\frac{\Lambda(\lambda_2,\mu_2)}{d_{1,i-1}^{\alpha_1} d_{2,j}^{\alpha_2}} \,T_A(T_{-1} \otimes I)QR_{i-1,j}f^{2} \nonumber\\
&+\frac{\Lambda(\lambda_3,\mu_3)}{d_{1,i-1}^{\alpha_1} d_{2,j-1}^{\alpha_2}}\,T_A( T_{-1}\otimes T_{-1})QR_{i-1,j-1}f^{3}+
\frac{\Lambda(\lambda_4,\mu_4)}{d_{1,i}^{\alpha_1} d_{2,j-1}^{\alpha_2}}\,T_A (I\otimes T_{-1})QR_{i,j-1}f^{4}.
\end{align}
}
Let us define 16 by 16 matrices $A_1$, $A_2$, $A_3$, $A_4$:
{\small
\begin{align*}
&  A_1=
\left[
  \begin{array}{c}
   \Lambda(\lambda_1,\mu_1)\,T_{A}QR_{i,j} \\
   \frac{\Lambda(\lambda_1,\mu_1)}{d_{1,i}^{}}\,T_{A}QR_{i,j} \\
   \frac{\Lambda(\lambda_1,\mu_1)}{d_{2,j}^{}}\,T_{A}QR_{i,j} \\
  \frac{\Lambda(\lambda_1,\mu_1)}{d_{1,i}^{} d_{2,j}^{}}\,T_{A}QR_{i,j} \\
  \end{array}
\right],\quad
  A_2=
\left[
  \begin{array}{c}
   \Lambda(\lambda_2,\mu_2)\,T_{A}QR_{i,j} \\
   \frac{\Lambda(\lambda_2,\mu_2)}{d_{1,i-1}^{}}\,T_{A}(T_{-1}\otimes I)QR_{i-1,j} \\
   \frac{\Lambda(\lambda_2,\mu_2)}{d_{2,j}^{}}\,T_{A}(T_{-1}\otimes I)QR_{i-1,j} \\
  \frac{\Lambda(\lambda_2,\mu_2)}{d_{1,i-1}^{} d_{2,j}^{}}\,T_{A}(T_{-1}\otimes I)QR_{i-1,j}\\
  \end{array}
\right],\\
&A_3=
\left[
  \begin{array}{c}
    \Lambda(\lambda_3,\mu_3)\,T_A( T_{-1}\otimes T_{-1})QR_{i-1,j-1} \\
    \frac{\Lambda(\lambda_3,\mu_3)}{d_{1,i-1}^{}}\,T_A( T_{-1}\otimes T_{-1})QR_{i-1,j-1} \\
    \frac{\Lambda(\lambda_3,\mu_3)}{d_{2,j-1}^{}}\,T_A( T_{-1}\otimes T_{-1})QR_{i-1,j-1} \\
   \frac{\Lambda(\lambda_3,\mu_3)}{d_{1,i-1}^{} d_{2,j-1}^{}}\,T_A( T_{-1}\otimes T_{-1})QR_{i-1,j-1} \\
  \end{array}
\right],\quad
A_4=
\left[
  \begin{array}{c}
   \Lambda(\lambda_4,\mu_4)\,T_A (I\otimes T_{-1})QR_{i,j-1} \\
    \frac{\Lambda(\lambda_4,\mu_4)}{d_{1,i}}\,T_A (I\otimes T_{-1})QR_{i,j-1} \\
    \frac{\Lambda(\lambda_4,\mu_4)}{d_{2,j-1}}\,T_A (I\otimes T_{-1})QR_{i,j-1} \\
    \frac{\Lambda(\lambda_4,\mu_4)}{d_{1,i}^{\alpha_1} d_{2,j-1}^{\alpha_2}}\,T_A (I\otimes T_{-1})QR_{i,j-1} \\
  \end{array}
\right],
\end{align*}
}
and let us define 4 by 4 matrices $a^{i,j}_k$, $k=1,\ldots,9$:
\begin{align*}
  & t_1=1:4,\quad t_2=5:8,\quad t_3=9:12,\quad t_4=13:16,\\
  &a^{i,j}_1 =A_3(:,t_1),\quad a^{i,j}_2=A_3(:,t_2)+A_4(:,t_1),\quad
  a^{i,j}_3=A_4(:,t_2),\quad  a^{i,j}_4=A_3(:,t_4)+A_2(:,t_1),\\
  & a^{i,j}_5=A_1(:,t_1)+A_2(:,t_2)+A_3(:,t_3)+A_4(:,t_4),\\
  &a^{i,j}_6=A_1(:,t_2)+A_4(:,t_3),\quad a^{i,j}_7=A_2(:,t_4),\quad a^{i,j}_8=A_1(:,t_4)+A_2(:,t_3),\quad a^{i,j}_9=A_1(:,t_3).
\end{align*}
Here we use the Matlab notation to express a sub matrix of a given matrix, for instance, for 4 by 16 matrix $A$ and the index $t_1=1:4$, we denote by $A(:,t_1)$ the sub matrix $[A_{i,j}]$, $i=1,\ldots,4$, $j=1,\ldots,4$ of the matrix $A$.

Let us label the multi-moment vectors at nearest nine grid points:
\begin{equation*}
  \mathbf{f}_{1}=\mathbf{f}_{i-1,j-1},\mathbf{f}_{2}=\mathbf{f}_{i,j-1},\mathbf{f}_{3}=\mathbf{f}_{i+1,j-1},
  \mathbf{f}_{4}=\mathbf{f}_{i-1,j},\mathbf{f}_{5}=\mathbf{f}_{i,j},\mathbf{f}_{6}=\mathbf{f}_{i+1,j},
  \mathbf{f}_{7}=\mathbf{f}_{i-1,j+1},\mathbf{f}_{8}=\mathbf{f}_{i,j+1},\mathbf{f}_{9}=\mathbf{f}_{i+1,j+1},
\end{equation*}
Then the left hand side of \eqref{formula1} for $(\alpha_1,\alpha_2)=(0,0),(1,0),(0,1),(1,1)$
are expressed using the nine vectors and the matrices $a^{i,j}_k$:
{\small
\begin{equation}\label{a}
\left[
  \begin{array}{c}
  \displaystyle  \sum_{k=1}^4 L((1+\xi\cdot \nabla)F_{k}(x_{i,j}+\cdot) | B_k) \\
\displaystyle    \sum_{k=1}^4 L((1+\xi\cdot \nabla)\partial_{\xi_1}F_{k}(x_{i,j}+\cdot) | B_k) \\
\displaystyle    \sum_{k=1}^4  L((1+\xi\cdot \nabla)\partial_{\xi_2}F_{k}(x_{i,j}+\cdot) | B_k) \\
 \displaystyle   \sum_{k=1}^4  L((1+\xi\cdot \nabla)\partial^2_{\xi_,\xi_2}F_{k}(x_{i,j}+\cdot) | B_k)\\
  \end{array}
\right]=
\sum_{k=1}^9a^{i,j}_k \mathbf{f}_k.
\end{equation}
}

As for $A=\partial_{\xi_1}\partial^{\alpha_1}_{\xi_1} \partial^{\alpha_2}_{\xi_2}$ and $A=\partial_{\xi_2}\partial^{\alpha_1}_{\xi_1} \partial^{\alpha_2}_{\xi_2}$, we obtain
{\small
\begin{align*}
&c\Delta t\sum_{k=1}^4 L(A F_{k}(x_{i,j}+\cdot) | B_k)\\
&=
\frac{\lambda_1\Lambda(\lambda_1,\mu_1)}{d_{1,i}^{\alpha_1} d_{2,j}^{\alpha_2}}\,T_{A}QR_{i,j} f^{1}
+
\frac{\lambda_2\Lambda(\lambda_2,\mu_2)}{d_{1,i-1}^{\alpha_1} d_{2,j}^{\alpha_2}} T_A(T_{-1} \otimes I)QR_{i-1,j}f^{2}\nonumber \\
&+\frac{\lambda_3\Lambda(\lambda_3,\mu_3)}{d_{1,i-1}^{\alpha_1} d_{2,j-1}^{\alpha_2}}\,T_A( T_{-1}\otimes T_{-1})QR_{i-1,j-1}f^{3}
+
\frac{\lambda_4\Lambda(\lambda_4,\mu_4)}{d_{1,i}^{\alpha_1} d_{2,j-1}^{\alpha_2}}\,T_A (I\otimes T_{-1})QR_{i,j-1}f^{4},\nonumber
\end{align*}
}
and
{\small
\begin{align*}
&c\Delta t\sum_{k=1}^4 L(A F_{k}(x_{i,j}+\cdot) | B_k)\\
&=
\frac{\mu_1\Lambda(\lambda_1,\mu_1)}{d_{1,i}^{\alpha_1} d_{2,j}^{\alpha_2}}\,T_{A}QR_{i,j} f^{1}
+
\frac{\mu_2\Lambda(\lambda_2,\mu_2)}{d_{1,i-1}^{\alpha_1} d_{2,j}^{\alpha_2}} T_A(T_{-1} \otimes I)QR_{i-1,j}f^{2}\nonumber \\
&+\frac{\mu_3\Lambda(\lambda_3,\mu_3)}{d_{1,i-1}^{\alpha_1} d_{2,j-1}^{\alpha_2}}\,T_A( T_{-1}\otimes T_{-1})QR_{i-1,j-1}f^{3}
+
\frac{\mu_4\Lambda(\lambda_4,\mu_4)}{d_{1,i}^{\alpha_1} d_{2,j-1}^{\alpha_2}}\,T_A (I\otimes T_{-1})QR_{i,j-1}f^{4}.\nonumber
\end{align*}
}
Similarly, one has
{\small
\begin{equation}\label{bc}
c\Delta t\left[
  \begin{array}{c}
  \displaystyle  \sum_{k=1}^4 L(\partial_{\xi_1}F_{k}(x_{i,j}+\cdot) | B_k) \\
\displaystyle    \sum_{k=1}^4 L(\partial_{\xi_1}\partial_{\xi_1}F_{k}(x_{i,j}+\cdot) | B_k) \\
\displaystyle    \sum_{k=1}^4 L(\partial_{\xi_1}\partial_{\xi_2}F_{k}(x_{i,j}+\cdot) | B_k) \\
 \displaystyle   \sum_{k=1}^4 L(\partial_{\xi_1}\partial^2_{\xi_,\xi_2}F_{k}(x_{i,j}+\cdot) | B_k)\\
  \end{array}
\right]=
\sum_{k=1}^9b^{i,j}_k \mathbf{f}_k,\quad
c\Delta t\left[
  \begin{array}{c}
  \displaystyle  \sum_{k=1}^4 L(\partial_{\xi_2}F_{k}(x_{i,j}+\cdot) | B_k) \\
\displaystyle    \sum_{k=1}^4 L(\partial_{\xi_2}\partial_{\xi_1}F_{k}(x_{i,j}+\cdot) | B_k) \\
\displaystyle    \sum_{k=1}^4  L(\partial_{\xi_2}\partial_{\xi_2}F_{k}(x_{i,j}+\cdot) | B_k) \\
 \displaystyle   \sum_{k=1}^4  L(\partial_{\xi_2}\partial^2_{\xi_,\xi_2}F_{k}(x_{i,j}+\cdot) | B_k)\\
  \end{array}
\right]=
\sum_{k=1}^9c^{i,j}_k \mathbf{f}_k,\quad
\end{equation}
}
with 4 by 4 matrices $b^{i,j}_k$, $c^{i,j}_k$, $k=1,\ldots,9$ which are defined in the same manner.
\subsection{The multi-moment scheme}
Suppose the numerical approximations to the exact solutions
$H^{n}_z(x_{i,j})$ and its first derivatives
$\partial_{x_1}H^{n}_z(x_{i,j})$ and
$\partial_{x_2}H^{n}_z(x_{i,j})$, and the second derivative
$\partial^2_{x_1x_2}H^{n}_z(x_{i,j})$ are known at all grid points
$x_{i,j}=(x_i,y_j)$ at time step $t_n$, which we denote by
\[
h^n_{ij},\quad \partial_xh^n_{ij},\quad \partial_yh^{n}_{ij},\quad
\partial^2_{x,y}h^{n}_{ij}.
\]
Similarly, for numerical approximations to
$E_x$ and $E_y$ and the derivatives, we use symbols $e^n_{x,i,j}$
and $e^n_{y,i,j}$.
Let $\mathbf{h}^n_{i,j}$,
${\mathbf{e_x}}^n_{i,j}$ and ${\mathbf{e_y}}^n_{i,j}$ denote the
multi moment vectors at grid; 
\begin{equation*} \label{alg}
\begin{array}{ll}
\mathbf{h}^n_{ij}&=[h^n_{ij},\partial_{x_1}h^n_{ij},\partial_{x_2}h^n_{ij},
\partial^2_{x_1x_2}h^n_{ij}]^\top,\\
{\mathbf{e_x}}^n_{ij}&=[{e_x}^n_{ij},\partial_{x_1}{e_x}^n_{ij},\partial_{x_2}{e_x}^n_{ij},
\partial^2_{x_1x_2}{e_x}^n_{ij}]^\top,\\
{\mathbf{e_y}}^n_{ij}&=[{e_y}^n_{ij},\partial_{x_1}{e_y}^n_{ij},
\partial_{x_2}{e_y}^n_{ij},\partial^2_{x_1x_2}{e_y}^n_{ij}]^\top.
\end{array} \end{equation*}
We number the numerical solutions in the way as above. Based on the exact integration formula \eqref{Mult}, and \eqref{a},\eqref{bc}, we arrive at the multi-moment scheme for the Maxwell's equations:
\begin{equation} \label{alg0}
\begin{array}{lll}
\mathbf{h}^{n+1}_{5}=
&\sum_{k=1}^9a^{i,j}_k \mathbf{h}^n_{k} +
\frac{1}{c\mu}c^{i,j}_k\mathbf{e_x}^{n}_{k}
&-\frac{1}{c\mu}b^{i,j}_k\mathbf{e_y}^{n}_{k},\\ \\
\mathbf{e_x}^{n+1}_{5}=
&\sum_{k=1}^9\frac{1}{c\epsilon}c^{i,j}_k \mathbf{h}^n_{k} +
a^{i,j}_k\mathbf{e_x}^{n}_{k},&\\ \\
\mathbf{e_y}^{n+1}_{5}=
&\sum_{k=1}^9-\frac{1}{c\epsilon}b^{i,j}_k \mathbf{h}^n_{k}&+a^{i,j}_k\mathbf{e_y}^{n}_{k},
\end{array}
\end{equation}
Thus, the
method uses 9 nearest neighbor points for 12 components (4 moments)
for $H_z$, $E_x$, and $E_y$. Each $4\times 36$ matrices $[a^{i,j}_1,\ldots,a^{i,j}_9]$,
 $[b^{i,j}_1,\ldots,b^{i,j}_9]$ and $[c^{i,j}_1,\ldots,c^{i,j}_9]$ has 100 nonzero entries.
Hence the total cost for the update \eqref{alg0} amounts to 700.

We would like to emphasis that the multi-moment scheme provides with the first and the second derivatives as we as the function value at each grid. When displaying the numerical solution, one can construct the bi-cubic polynomial and can evaluate any spacial point with using the interpolation.

\section{Stability}\label{sec:stability}
We analyze the stability of the multi-moment scheme. For simplicity, we assume that the grid length is uniform, i.e.,
\begin{equation*}
  \Delta x:= x_{i}-x_{i-1}=y_{j}-y_{j-1}
\end{equation*}
for all $i,j$.
In this case, the 4 by 4 matrices $a^{i,j}_k$, $b^{i,j}_k$ and $c^{i,j}_k$ in \eqref{alg0} remain the same throughout of $(i,j)$, and thus we omit the subscript.
Suppose that $\mathbf{f}=\{\mathbf{f}_{\ell,j}\}_{\ell,j}$ is $N$ periodic
with respect to $\ell,j$, i.e.,
\[
\mathbf{f}_{0,j}=\mathbf{f}_{N,j}, \quad \mathbf{f}_{\ell,0}=\mathbf{f}_{\ell,N},
\]
for all $0\le \ell,\; j\le M$. Let us consider the discrete Fourier transform:
\begin{equation*}
\widehat{\mathbf{f}}[k] :=\sum_{i=0}^N\sum^N_{j=0} \mathbf{f}_{i,j}
e^{-2\pi i k\cdot x_{i,j}}
=\sum_{i=0}^N\sum^N_{j=0}\left[
                                 f_{i,j},
                                 \partial_{x_1}f_{i,j},
                                 \partial_{x_2}f_{i,j},
                                 \partial^2_{x_1x_2}f_{i,j}
 \right]^\top
e^{-2\pi i k\cdot x_{i,j}},
\end{equation*}
where $k\cdot x_{i,j}=(k_1,k_2) \cdot (x_i,y_j) = k_1 x_i +
k_2 y_j$. The discrete Fourier transform of the sequence
$\mathbf{F}_{i,j}:=
a_1 \mathbf{f}_{i-1,j-1}+ a_2 \mathbf{f}_{i,j-1}+a_3 \mathbf{f}_{i+1,j-1}+
a_4 \mathbf{f}_{i-1,j}+a_5 \mathbf{f}_{i,j}+a_6 \mathbf{f}_{i+1,j}+
a_7 \mathbf{f}_{i-1,j+1}+a_8 \mathbf{f}_{i,j+1}+a_9 \mathbf{f}_{i+1,j+1}
$ becomes:
\begin{equation}\label{fourier}
\sum_{\ell=0}^N\sum^N_{j=0} \mathbf{F}_{i,j}
e^{-2\pi i k\cdot x_{\ell,j}} =g(\{a_k\},2\pi k_1\Delta x,2\pi k_2\Delta
x) \widehat{\mathbf{f}}[k],
\end{equation}
where $g(\{a_k\},2\pi k_1\Delta x,2\pi k_2\Delta x)$ is $4 \times 4$ matrix depending on $\lambda$, $k_1\Delta x$ and $k_2\Delta x$:
\begin{align*}
& g(\{a_k\},2\pi k_1\Delta x,2\pi k_2\Delta x)
 =e^{2\pi i(-k_1-k_2)\Delta x} a_1
  +e^{-2\pi i k_2\Delta x} a_2
  +e^{2\pi i(k_1-k_2)\Delta x} a_3  \\
&+e^{-2\pi i k_1\Delta x}a_4
+ a_5
+e^{2\pi i k_1 \Delta x}a_6
+e^{2\pi i(-k_1+k_2)\Delta x}a_7
+e^{2\pi i k_2\Delta x} a_8
+ e^{2\pi i(k_1+k_2)\Delta x}a_9.
\end{align*}
Now let us assume $(H_z,E_x,E_y)$ is periodic. Let
$\mathbf{h}^n=\{\mathbf{h}^{n}_{i,j}\}_{i,j}$, $\mathbf{e_x}^n=\{\mathbf{e_x}^{n}_{i,j}\}_{i,j}$ and  $\mathbf{e_y}^n=\{\mathbf{e_y}^{n}_{i,j}\}_{i,j}$.
From \eqref{alg0} and \eqref{fourier}, the amplification factor $\Phi(\lambda,2\pi k_1\Delta x,2\pi k_2\Delta x)$ (12 by 12 matrix) of the proposed scheme is given by
\begin{equation}\label{Phi}
\Phi(\lambda,\theta_1,\theta_2) =
\left[
\begin{array}{ccc}
g(\{a_k\},\theta_1,\theta_2) & \frac{1}{c\mu}g(\{c_k\},\theta_1,\theta_2) & -\frac{1}{c\mu}g(\{b_k\},\theta_1,\theta_2) \\
 \frac{1}{c\varepsilon}g(\{c_k\},\theta_1,\theta_2) & g(\{a_k\},\theta_1,\theta_2) & \mathbb{O} \\
 -\frac{1}{c\varepsilon}g(\{b_k\},\theta_1,\theta_2) & \mathbb{O} & g(\{a_k\},\theta_1,\theta_2)
 \end{array}
 \right],
\end{equation}
where $\theta_1 = 2\pi k_1 \Delta x$, $\theta_2 = 2\pi k_2 \Delta x$, i.e., $\Phi(\lambda,2\pi k_1\Delta x,2\pi k_2\Delta x)$ maps
$
(\widehat{\mathbf{h}^{n}}, \widehat{\mathbf{e_x}^{n}},\widehat{\mathbf{e_y}^{n}})
$ to the next step
$
(\widehat{\mathbf{h}^{n+1}}, \widehat{\mathbf{e_x}^{n+1}},\widehat{\mathbf{e_y}^{n+1}})
$.

\begin{figure}[h]
\centering
\includegraphics[width=5cm]{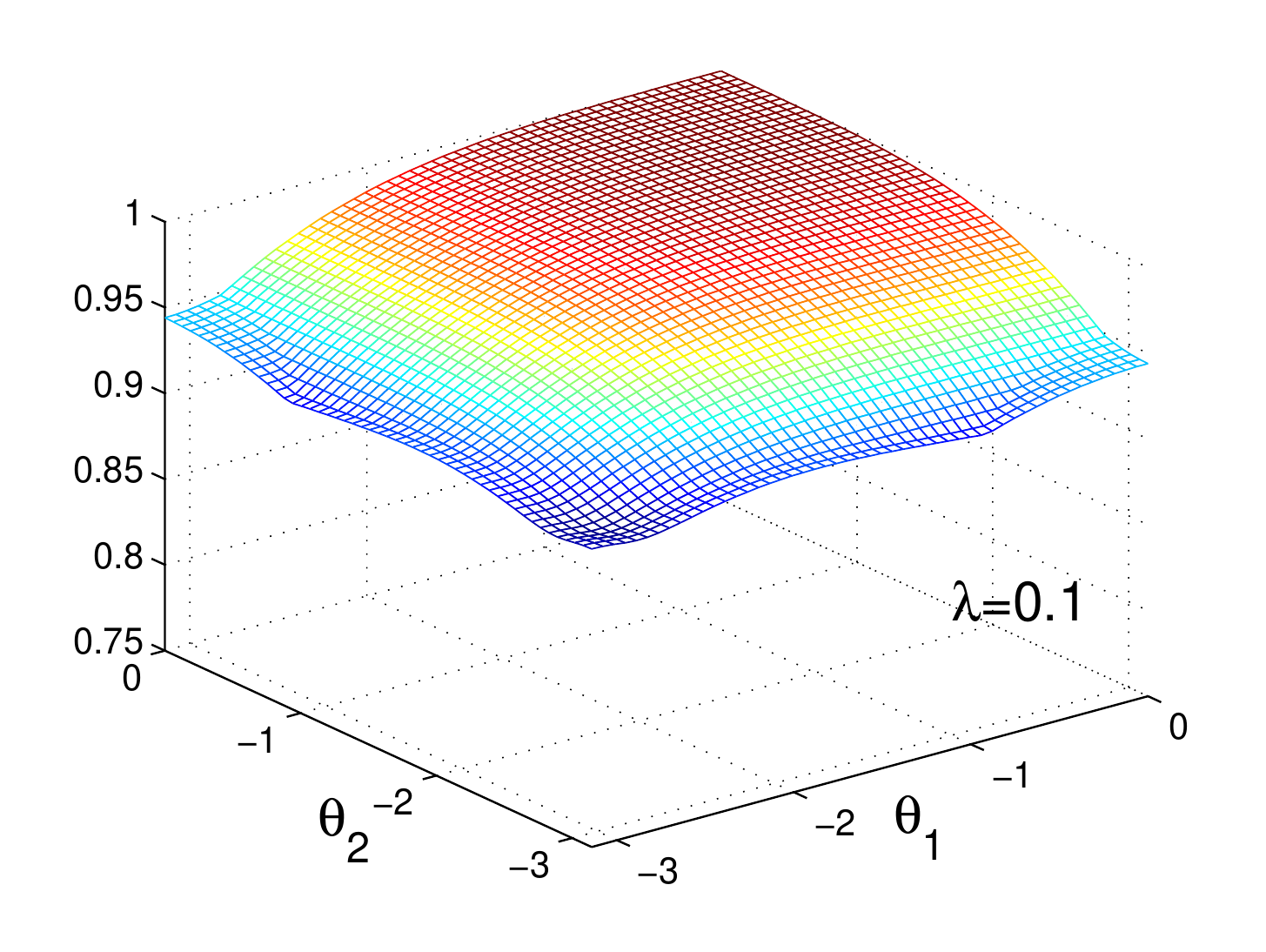}
\includegraphics[width=5cm]{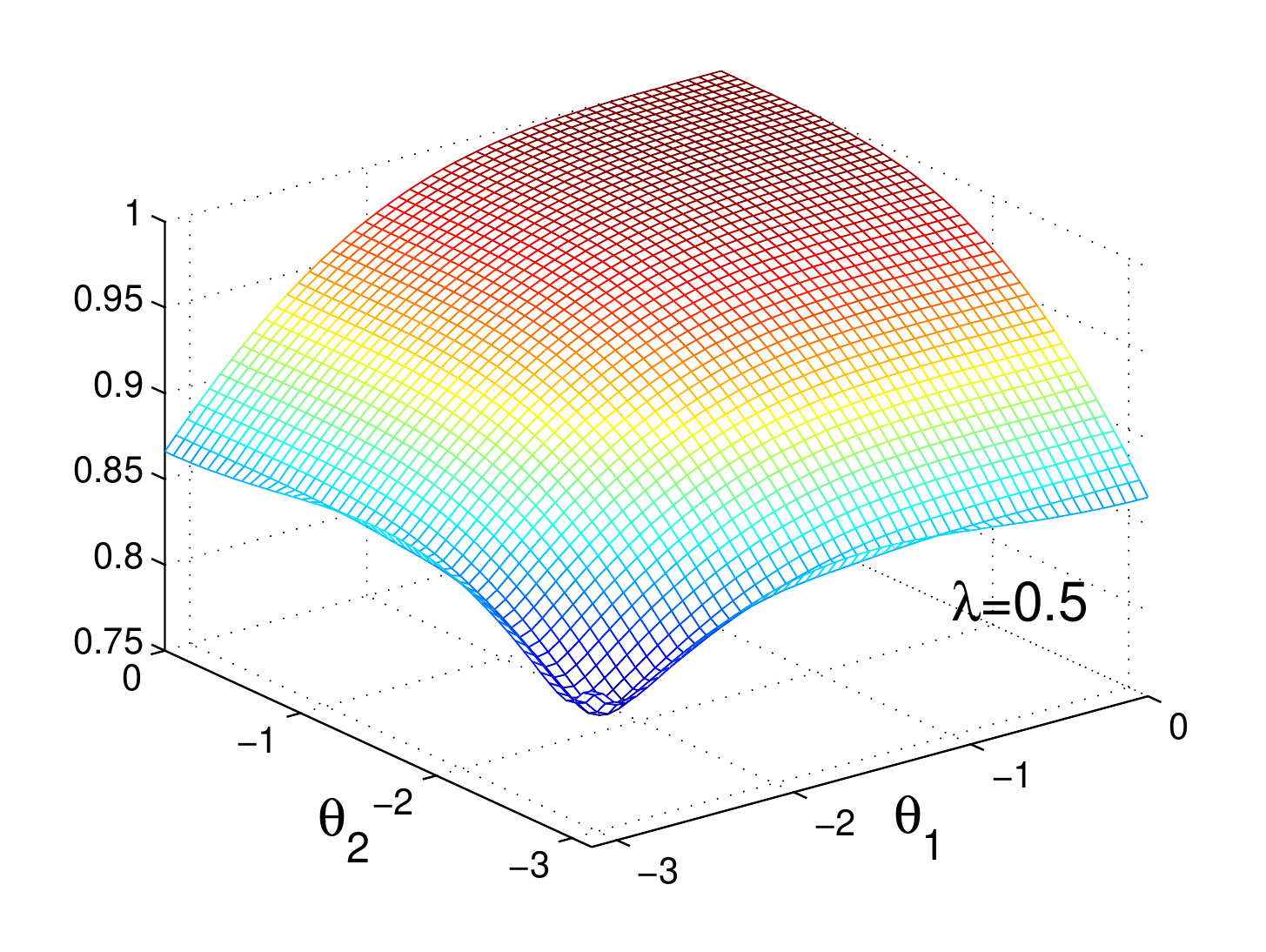}
\includegraphics[width=5cm]{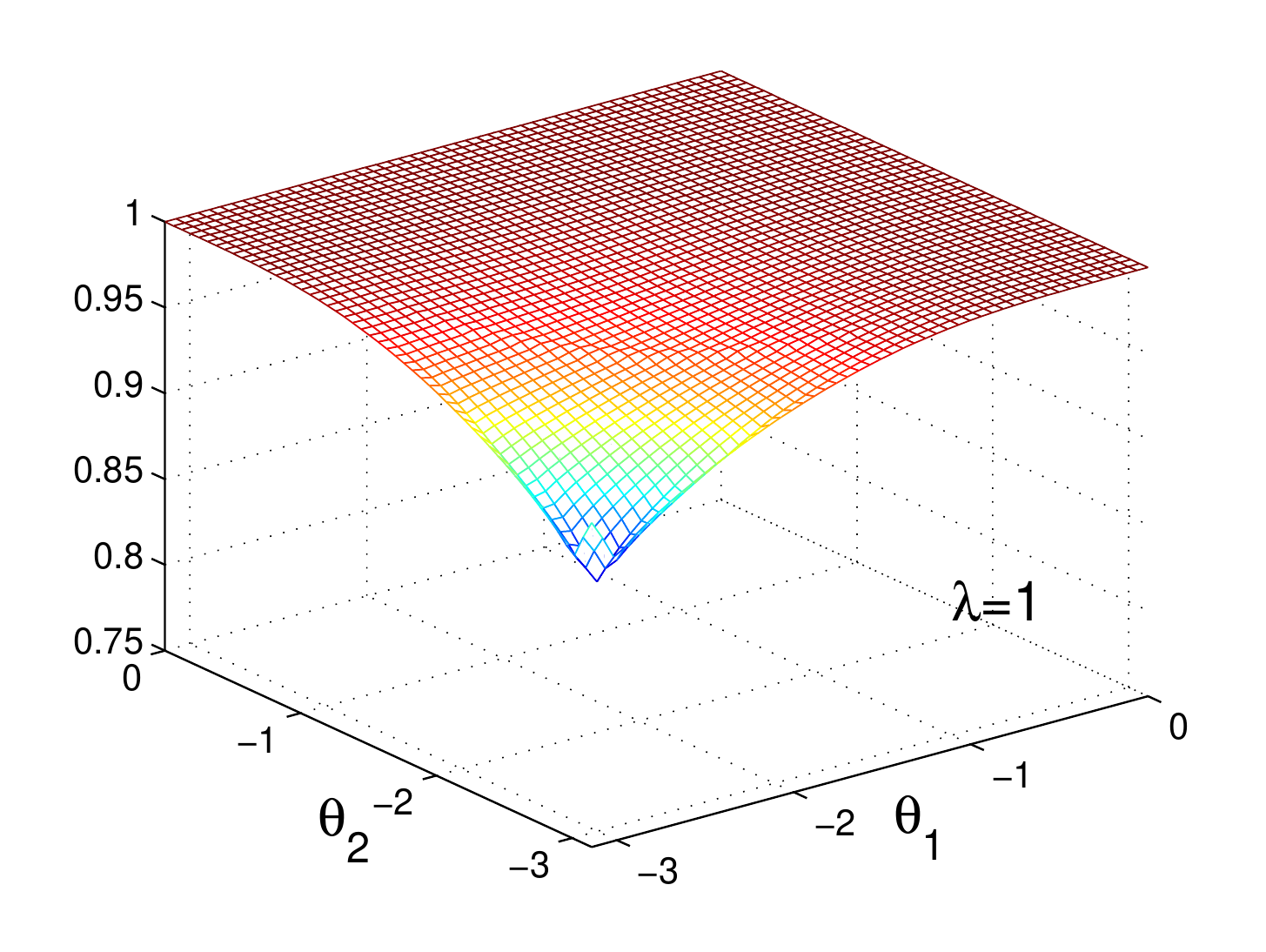}
\caption{The maximum absolute value of the eigenvalues of the amplification matrix \eqref{Phi} against $(\theta_1,\theta_2)$. $\lambda=0.2$ (left), $\lambda=0.5$ (center), $\lambda=1$, (right). }\label{fig:amp1}
\end{figure}

Let $\rho(\lambda,\theta_1,\theta_2)$ denote the set of eigenvalues of $\Phi(\lambda,\theta_1,\theta_2)$ ( 9 eigenvalues).
Figure \ref{fig:amp1} shows the maximum absolute value of the eigenvalues, $\max\{|\rho(\lambda,\theta_1,\theta_2)|\}$, against $[\theta_1,\theta_2]\in [-\pi,0]\times [-\pi,0]$ for   $\lambda=0.2$, $\lambda=0.5$ and $\lambda=1$, respectively. Numerically we find that all eigenvalues have the magnitude equal to or less than 1 for arbitrary $(\theta_1,\theta_2)$. The magnitude is close to 1 in a wide range of $[\theta_1,\theta_2]$, which indicates that the numerical scheme is less dissipative.

\section{Numerical test}\label{sec:numericaltest}
In this section, we show the numerical performance of the multi-moment scheme through some numerical tests.
\\
{\bf Example 1. Plane waves.}\\
In this example, we compute the numerical solutions for plane waves.
We compare the numerical solutions with those produced by the fourth order in time and space FDTD (Yee's scheme).
The Yee's scheme computes E field and H field at different time level, and thus one must provide the exact initial condition at time $t=-\frac{\Delta t}{2}$ for E field and $t=0$ for H field to obtain an accurate numerical solution. The plane wave solution is suitable to avoid the issue with the initial condition since the exact solution is easily obtained.

Let $f_\sigma(x,y) =\exp(- \frac{x^2}{\sigma^2} )$ for $\sigma>0$. Let $f_\sigma$ also denote its periodic extension to $x$ direction with periodicity $L$, i.e.,
$f_\sigma(x + L,y) = f_\sigma(x,y)$ for all $(x,y) \in\R^2$.
We rotate the function with the angle $\theta$ with respect to the origin $(0,0)$ to construct the one way propagating plane wave solution for Maxwell's equation:
\begin{equation*}
  H_z(t,x,y) =f_\sigma(x\cos\theta -y\sin\theta- t ) ,\quad
  E_x(t,x,y) = H_z(t,x,y)\sin\theta , \quad
  E_y(t,x,y) = H_z(t,x,y)\cos\theta.
\end{equation*}
If we set $L=\cos\theta $ with $\theta =\tan^{-1}m$ for $m\in \N$, they
are the periodic solution of the Maxwell's equation for $\epsilon=\mu=1$ in the domain $\Omega=[-\frac{1}{2},\frac{1}{2}]\times[-\frac{1}{2},\frac{1}{2}]$. In this numerical test we consider the solutions
{\small
\begin{align*}
  &H_z(t,x,y)=\sum_{m=0}^3f_\sigma(x\cos\theta_m -y\sin\theta_m- t ),\\
  &E_x(t,x,y) = \sum_{m=0}^3H_z(t,x,y)\sin\theta_m,\quad
  E_y(t,x,y) = \sum_{m=0}^3H_z(t,x,y)\cos\theta_m,
\end{align*}
}
where $\theta_m =\tan^{-1}m$.
We test the case $\sigma^{-2}=500, 1500$. Figure 1 shows the initial profile of $H_z(0,x,y)$, and the four plane waves at $t=0$, $f_\sigma(x\cos\theta_m -y\sin\theta_m)$ for $m=0,1,2,3,$ with $\sigma^{-2}=500$, in the domain. The arrow in each plot shows the direction of wave propagation.

We report the accuracy of the multi moment scheme.
The time step size is fixed to be $\Delta t= \lambda \Delta x$ for each mesh size $\Delta x = N^{-1}$, $N\in\{50, 100, 150, 200\}$. We report the numerical solutions for $\lambda = 1$ and $\lambda=\frac{1}{\sqrt{2}}$.
The numerical solutions at time $T=1$ are produced by the multi moment scheme and compared to the exact solution. The number of iteration is $N$ for $\lambda =1$, and $1.4 N$ for $\lambda = \frac{1}{\sqrt{2}}$ where the numerical solution approximates the solution at time $T= \frac{1.4}{\sqrt{2}} \sim 0.98995$.

 The initial value for $\mathbf{h}^0_{i,j}$ is provided by the exact solution:
 {\small
\begin{equation*}
  h^0_{i,j} = H_z(0,x_i,y_j),\; \partial_xh^0_{i,j} = \partial_xH_z(0,x_i,y_j),\;
  \partial_yh^0_{i,j} = \partial_yH_z(0,x_i,y_j),\; \partial^2_{x,y}h^0_{i,j} = \partial^2_{x,y}H_z(0,x_i,y_j).
\end{equation*}
}
Similarly, $\mathbf{e_x}^0_{i,j}$ and $\mathbf{e_y}^0_{i,j}$ are given by the exact solution. In a practical situation, the initial condition in function form may not be available, and we only have the function value at each grid point. In that case we employ finite difference of the initial grid function to provide the initial condition for derivatives.

For each mesh size $\frac{1}{N}$, the error in the numerical solutions is measured by $\ell^\infty$ norm:
{\small
\begin{equation*}\label{CIPError}
  \epsilon_1=\max_{{i,j}}\{|h^n_{i,j} - H_z(T,x_{i,j})|,|(e_x)^n_{i,j} - E_x(T,x_{i,j})|,|(e_y)^n_{i,j} - E_y(T,x_{i,j})|\}
\end{equation*}
}
The multi-moment scheme computes the first derivatives and the second mix derivative as numerical solutions. We report the relative error in the first derivatives:
{\small
\begin{align*}\label{CIPError1st}
\epsilon_2
=\max_{{i,j},\alpha}\{\frac{|\partial_\alpha h^n_{i,j} - \partial_\alpha  H_z(T,x_{i,j})|}{ | \partial_\alpha  H_z(T,x_{i,j}) | },\frac{|(\partial_\alpha e_x)^n_{i,j} - \partial_\alpha E_x(T,x_{i,j})|}{| \partial_\alpha E_x(T,x_{i,j})| },\frac{|(\partial_\alpha e_y)^n_{i,j} - \partial_\alpha E_y(T,x_{i,j})|}{|\partial_\alpha E_y(T,x_{i,j}) | }\},
\end{align*}
}
where $\alpha=x,y$.

For comparison, we also report the numerical solution produced by fourth-order in time and space FDTD (see \cite{Deveze+BeaulieuETAL:92}). In the FDTD numerical simulation, we employ the CFL number to be $\frac{1}{\sqrt{2}}$ with which the FDTD method provides the best performance, i.e., FDTD with CFL$=\frac{1}{\sqrt{2}}$ produces most accurate numerical solutions among the other CFL.
In this test, the initial values for $h^0_{i,j}$, ${e_x}^{-\frac{1}{2}}_{i,j}$ and ${e_y}^{-\frac{1}{2}}_{i,j}$ are given exactly:
\begin{equation*}
 h^0_{i,j} = H_z(0,x_i,y_j),\quad  (e_x)^{-\frac{1}{2}}_{i,j} = E_x(-\tfrac{\Delta t}{2},x_i,y_j),\quad (e_y)^{-\frac{1}{2}}_{i,j} = E_y(-\tfrac{\Delta t}{2},x_i,y_j).
\end{equation*}
When the analytic solutions are not available, one must solve the Maxwell's equations backward in time to provide an accurate initial condition for $E_x(-\frac{\Delta t}{2},x_i,y_j)$ and $E_y(-\frac{\Delta t}{2},x_i,y_j)$ to start the FDTD scheme, by employing another time marching method such as Runge-Kutta schemes.

We compute the error in the numerical solution:
{\small
\begin{align*}
&  \epsilon_1=\max_{{i,j}}\{|h^n_{i,j} - H_z(T,x_{i,j})|,|(e_x)^n_{i,j} - E_x(T,x_{i,j})|,|(e_y)^n_{i,j} - E_y(T,x_{i,j})|\},\\
&\epsilon_2=\max_{{i,j},\alpha\in{x,y}}\{\frac{|\partial_\alpha h^n_{i,j} - \partial_\alpha  H_z(T,x_{i,j})|}{ | \partial_\alpha  H_z(T,x_{i,j}) | },\frac{|(\partial_\alpha e_x)^n_{i,j} - \partial_\alpha E_x(T,x_{i,j})|}{| \partial_\alpha E_x(T,x_{i,j})| },\frac{|(\partial_\alpha e_y)^n_{i,j} - \partial_\alpha E_y(T,x_{i,j})|}{|\partial_\alpha E_y(T,x_{i,j}) | }\},
\end{align*}
}
where $\partial_\alpha h^n_{i,j}$, $ (\partial_\alpha e_x)^n_{i,j}$ and $(\partial_\alpha e_y)^n_{i,j}$ are computed by employing third order finite difference which is used in the fourth order FDTD scheme.

Figure \ref{fig:Error500} shows the error $\epsilon_1$ and $\epsilon_2$ of the solutions generated by the multi moment method with $\lambda=1$ and $\lambda=\frac{1}{\sqrt{2}}$, and the fourth order FDTD with $\lambda=\frac{1}{\sqrt{2}}$ for the initial profile with $\sigma^{-2}=500$ against the grid number $N$.
The order of accuracy of the numerical solutions for each method is shown in Table \ref{table:order500}.
And Figure \ref{fig:Error1500} and Table \ref{table:order1500} are numerical results when $\sigma^{-2}=1500$.

Let us consider the total computational cost in the multi-moment scheme to obtain the numerical solution at $T=1$. For each time step, $300 N^2$ operation is required to update $H_z$, and $200 N^2$ for each $E_x$ and $E_y$. Since the number of time integration in the multi-moment scheme with $\lambda=1$ is $\frac{1}{\Delta t}$, the total cost (operation count) amounts to $700 N^2\times N = 700 N^3$.

Let us focus on the numerical solutions by the multi-moment scheme and Yee's scheme when $N=50$.
Figure \ref{fig:Error500} shows that the error in the numerical solution produced by the multi-moment scheme with $\lambda=1$ is $10^{-2}$ while the error by FDTD is $10^{-1}$. So, the multi-moment scheme is 10 times more accurate than the FDTD with $\lambda=\frac{1}{\sqrt{2}}$ for this mesh size. To obtain the same accurate numerical solution by FDTD, we must take the half mesh size $1/(2N)=1/100$. For the fourth order Yee's scheme, the cost for the one step map is 51 and the number of time integration is $\frac{1}{\Delta t}=2\sqrt{2}N$, thus the total cost for the numerical solution at $T=1$ amounts to $2\sqrt{2} N \times 51 (2N)^2 \sim 577 N^3$.


\noindent
{\bf Example 2. Sharp profile.}\\
Next we solve \eqref{Max} with $\epsilon=\mu=1$. The initial condition is
\begin{equation*}
  H(0,x,y)=
  \left\{
  \begin{array}{l}
    1,\quad x\in D,\\
    0,\quad x\in D^{\complement}.
  \end{array}\right.,\quad
  E_x(0,x,y)=E_y(0,x,y)=0,
\end{equation*}
where $D=[0.25, 0.75]\times[0.25, 0.75]$. The mesh size is $0.01$.
The initial condition is approximated by bi-cubic polynomial with the first and second derivatives being 0. We do not use any other techniques to approximate the initial discontinuous profile. We use our algorithm \eqref{alg0} with CFL$=1$.
The initial condition and the numerical solution for $H$ at time $T=0.15$ and $T=0.25$ are displayed in Figure \ref{fig:square}. No oscillation is observed in the numerical solution. We also employed the fourth order FDTD with the same initial condition. Numerical oscillations were found near the sharp profile. \\
{\bf Example 3. Hidden resolution}\\
We solve \eqref{Max} with $\epsilon=\mu=1$ with the initial condition is
\begin{equation*}
  H(0,x,y)=\exp\left(-\frac{x^2}{1000}\right),\quad
  E_x(0,x,y)=E_y(0,x,y)=0.
\end{equation*}
The mesh size is $\frac{1}{40}$ and CFL$=1$, and so $\Delta t = \frac{1}{40}$. We compute the numerical solutions at $t=10\Delta t$. We denote the numerical solutions by $\{ {h^{10}_z}_{,i,j} \}_{i,j}$, $\{ {\partial_x h^{10}_z}_{,i,j} \}_{,i,j}$, $\{ {\partial_y h^{10}_z}_{,i,j} \}_{i,j}$,
$\{{\partial^2_{x,y}h^{10}_z}_{i,j} \}_{i,j}$.
When visualizing the numerical solutions, we usually construct the bi-linear interpolation in each cell using the numerical solution at the grid.
For instance, for the visualization of the solution $\{ {h^{10}_z}_{,i,j} \}_{i,j}$, we plot the bi-linear interpolation constructed by using the gird value $\{ {h^{10}_z}_{,i,j} \}_{i,j}$, and for visualizing the numerical solution $\{ {\partial_x h^{10}_z}_{,i,j} \}_{i,j}$, we plot the bi-linear interpolation constructed by the gird value $\{ {\partial_x h^{10}_z}_{,i,j} \}_{i,j}$. These two bi-linear interpolations are unrelated.
In the left column in Figure \ref{fig:interp}, we plot the piece wise bi-linear interpolation for $\{ {h^{10}_z}_{i,j} \}_{i,j}$ (top) and the one for $\{ {\partial_x h^{10}_z}_{i,j} \}_{i,j}$ (bottom). One can observe the spiky peaks and dips in the plots.

As have been mentioned repeatedly, the multi-moment scheme produces the derivatives as well as the function value at each grid. This is equivalent to state that the multi-moment scheme computes the bi-cubic polynomial in each cell as a numerical solution. So when plotting the numerical solution, we should use the computed bi-cubic interpolation instead of bi-linear interpolation.
 Let us construct the bi-cubic polynomial ${h_z}_{,i,j}(x,y)$ in each cell $[x_{i-1},x_i]\times[y_{j-1},y_j]$ using the numerical solutions, and let $h_z(x,y)$ denote the piece wise bi-cubic polynomial defined in the domain $[-0.5,-0.5]\times [-0.5,-0.5]$, i.e.,
 \begin{equation}\label{bicubicHz}
   h_z(x,y) = {h_z}_{,i,j}(x,y), \quad (x,y) \in [x_{i-1},x_i]\times[y_{j-1},y_j].
 \end{equation}
In the right column of Figure \ref{fig:interp}, we show the surface plot of the bi-cubic interpolations $h_z(x,y)$ (top) and $\partial_xh_z(x,y)$ (bottom).
The numerical solutions are depicted as smooth functions. \\
{\bf{Derivative free method}}\\
We have also implemented the method using the bi-linear interpolation at each cell $C_{i,j}$:
$$
F(x,y)=F(0,0)(1-x)(1-y)+F(\Delta x ,0)x(1-y)+F(0,\Delta x)(1-x)y
+F(\Delta x,\Delta x)xy
$$
In this way we obtain a derivative free nine point scheme:
{\small
\begin{equation*}
\left[\begin{array}{c} h^{n+1}_{i,j} \\
(e_x)^{n+1}_{i,j} \\
(e_y)^{n+1}_{i,j}\end{array}\right] =\left[\begin{array}{ccc}
L_1 &\frac{1}{c\mu}L_2 & -\frac{1}{c\mu}L_3 \\
\frac{1}{c\varepsilon}L_2 & L_1 & 0\\
-\frac{1}{c\varepsilon}L_3 & 0 & L_1\end{array}\right]
\left[\begin{array}{c} \{h^{n}_{i,j}\} \\
\{(e_x)^{n}_{i,j}\} \\
\{(e_y)^{n}_{i,j}\}\end{array}\right],
\end{equation*}
}
where
{\small
\begin{align*}
&L_1=\left[\frac{\lambda ^2}{2 \pi },
\frac{(\pi -2 \lambda ) \lambda }{2 \pi },
\frac{\lambda ^2}{2 \pi } ,
\frac{(\pi -2 \lambda ) \lambda }{2 \pi } ,
 1-2 \lambda +\frac{2 \lambda ^2}{\pi } ,
 \frac{(\pi -2 \lambda ) \lambda }{2 \pi },
 \frac{\lambda ^2}{2 \pi } ,
 \frac{(\pi -2 \lambda ) \lambda }{2 \pi } ,
 \frac{\lambda ^2}{2 \pi } \right],\\
&L_2=\left[ -\frac{\lambda }{8} , -\frac{2-\lambda}{4} ,-\frac{\lambda }{8} ,
 0 , 0 , 0 ,\frac{\lambda}{8},\frac{2-\lambda}{4},\frac{\lambda }{8}\right], \quad L_3=\left[ -\frac{\lambda}{8} , 0 ,\frac{\lambda}{8} ,-\frac{2-\lambda}{4},0 , \frac{2-\lambda }{4} ,
 -\frac{\lambda }{8},0 ,\frac{\lambda }{8}\right].
\end{align*}
}
This method is also stable with $CFL=1$ but is second order accurate.
Our numerical tests show that if we let $\lambda=1$, (CFL$=1$) then
there is no significant dissipation but $\lambda=0.5$ it has 30\%
dissipation at $T=1$ with speed one. For oblique plane waves there
is no significant phase error with $CFL=1$. An advantage of this method is
that it is simple to be implemented and to be extended to the three dimension case.

\section{Conclusion}
We developed a numerical method for solving time-domain Maxwell's
equation. It is fully explicit space and
time integration method with higher order accuracy and CFL number
being one. The bi-cubic interpolation is used for the solution
profile to attain the resolution. It preserves sharp profiles very accurately
without any smearing and distortion due to
the exact time integration and high resolution approximation.
The stability of the method were analyzed, and the nearly forth order accuracy were observed.

\bibliographystyle{siam}

\begin{figure}[h]
\centering
\includegraphics[scale=0.2]{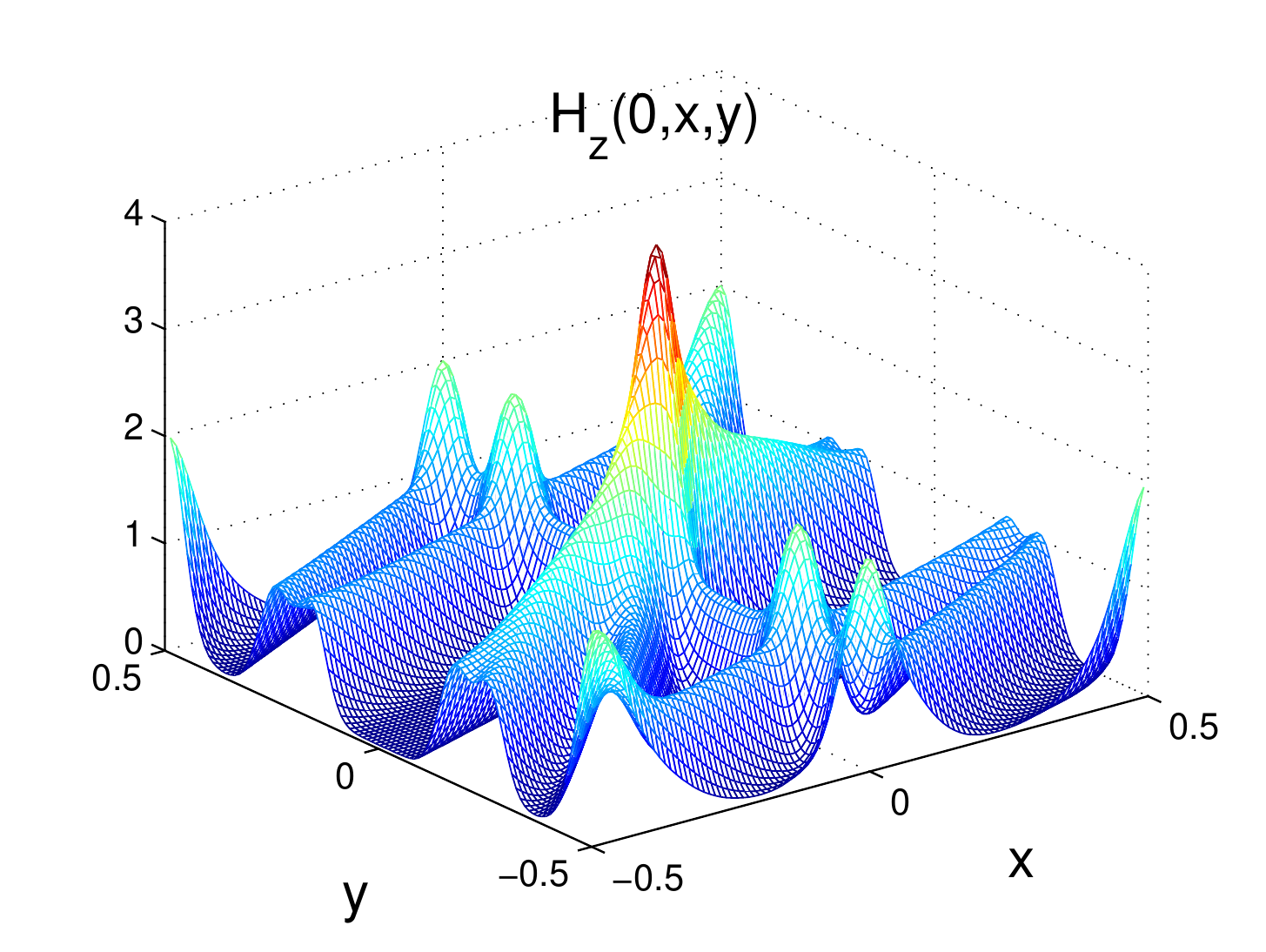}
\includegraphics[scale=0.2]{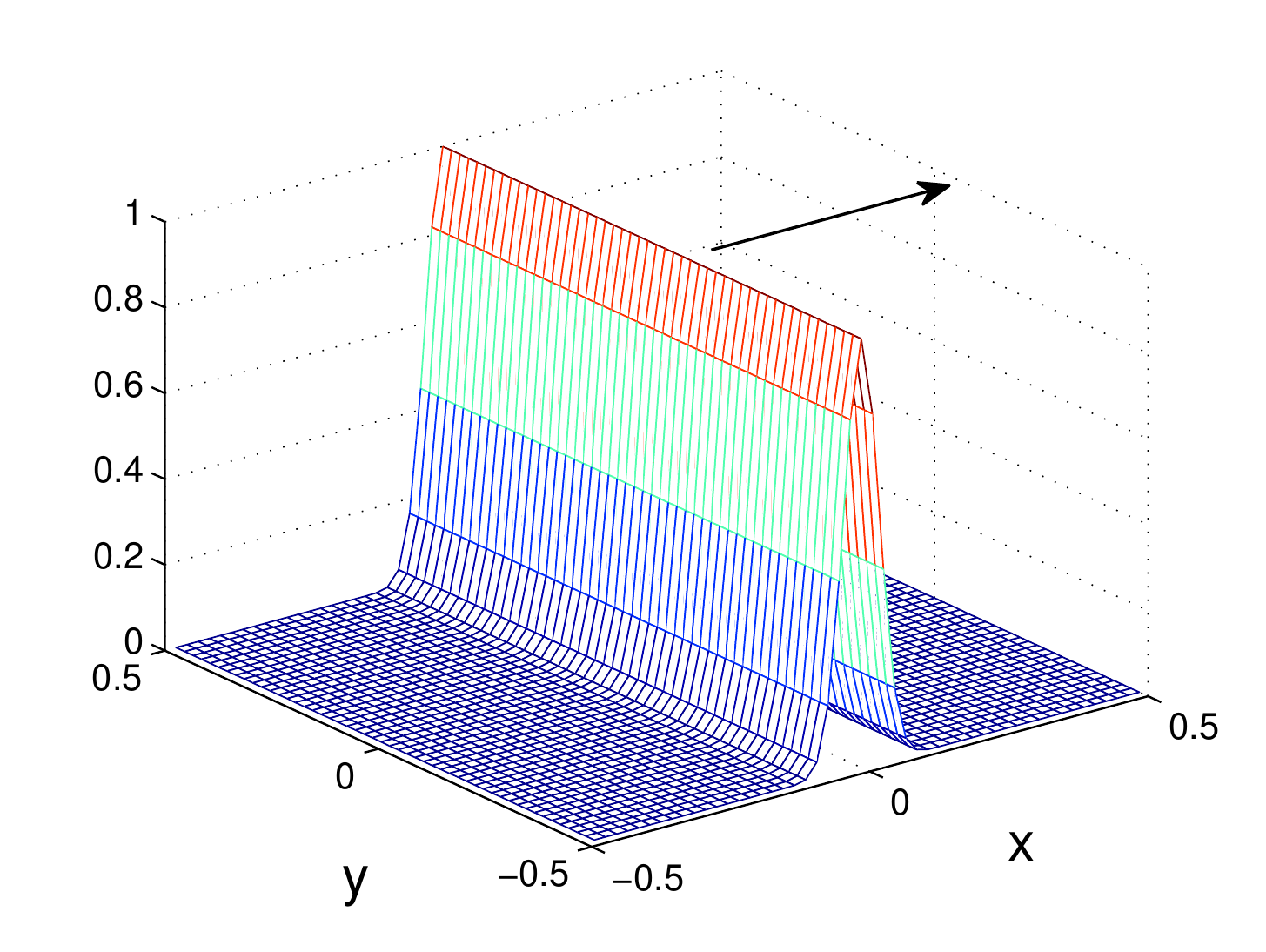}
\includegraphics[scale=0.2]{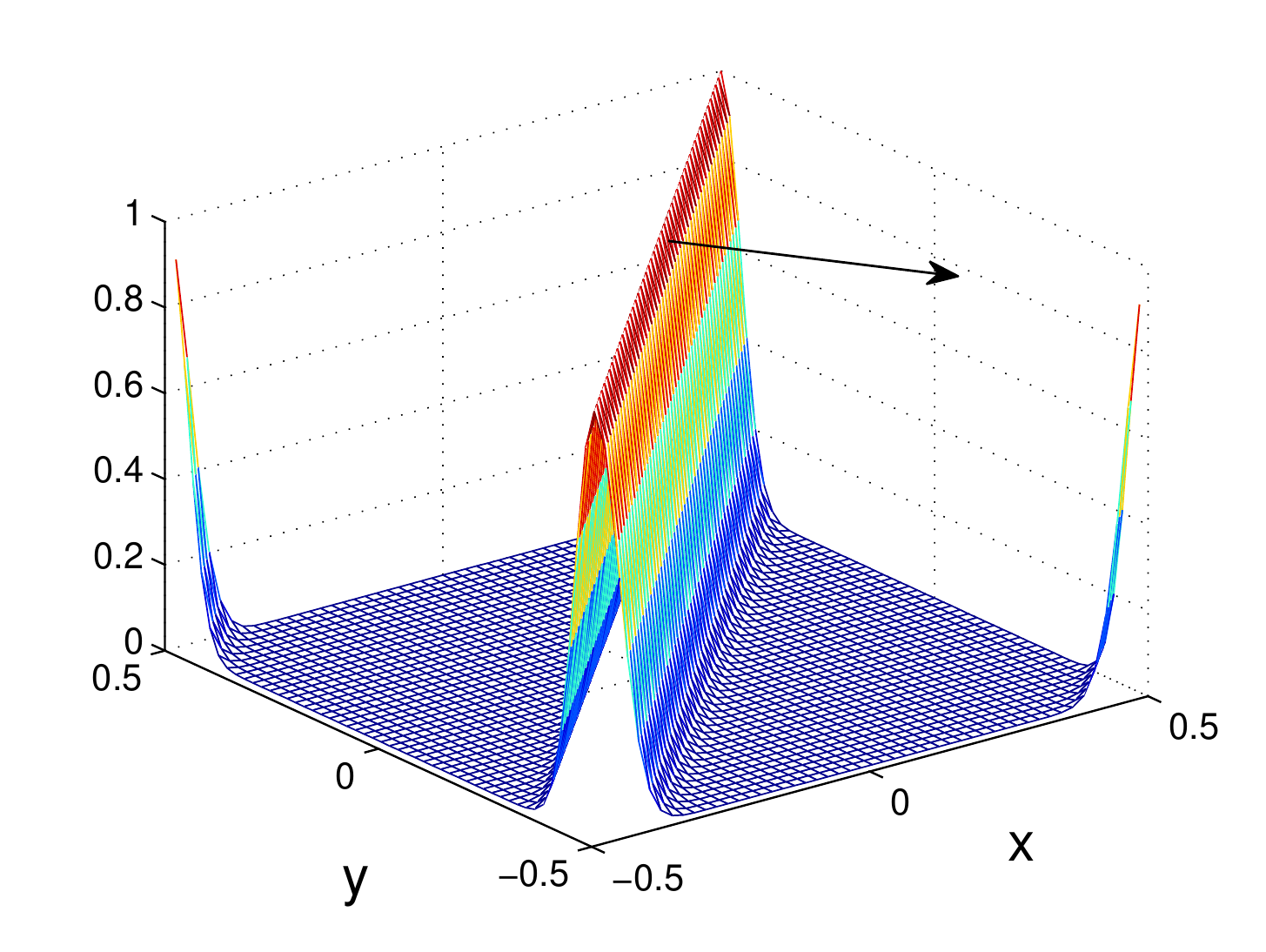}
\includegraphics[scale=0.2]{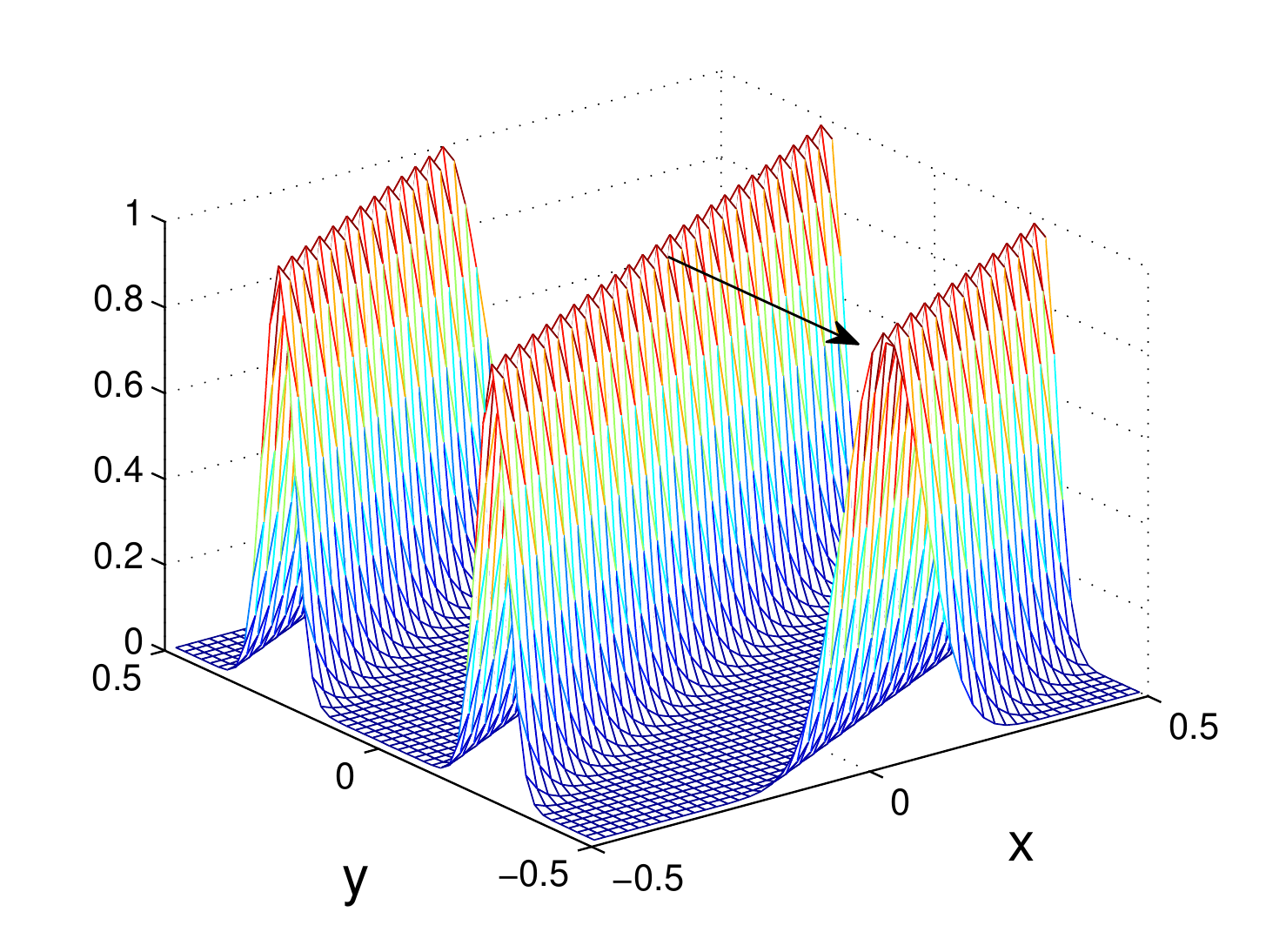}
\includegraphics[scale=0.2]{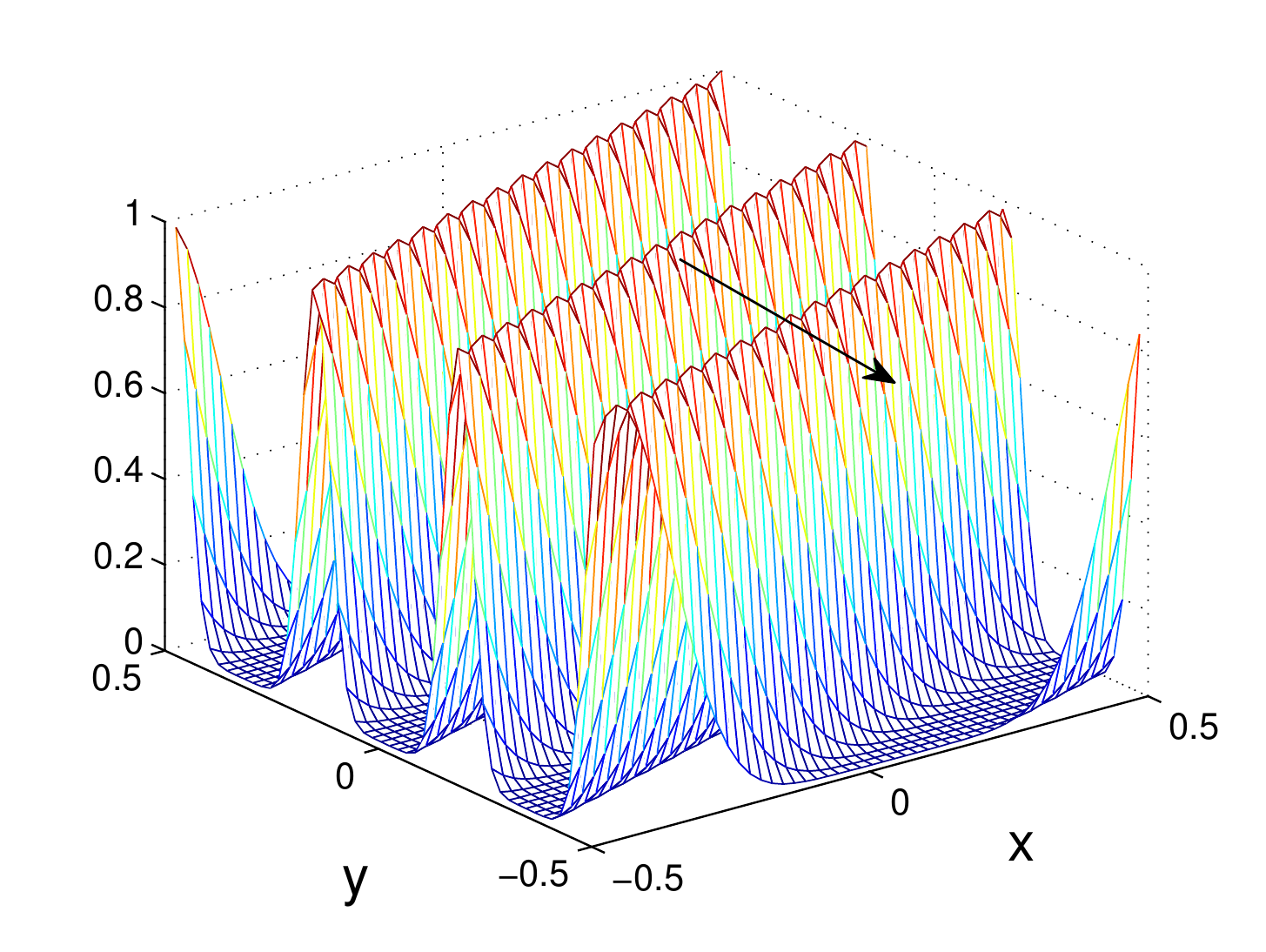}
\caption{Initial profiles $H_z(0,x,y)=\sum_{m=0}^3f_\sigma(x\cos\theta_m-y\sin\theta_m)$ (left) and $f_\sigma(x\cos\theta_m-y\sin\theta_m)$ for $m=0,1,2,3$. The inverse of the variance $\sigma^{-2}$ is $500$. The arrow in each plot shows the direction of wave propagation.}\label{fig:initial}
\end{figure}

\begin{figure}[h]
\centering
\includegraphics[width=7.5cm]{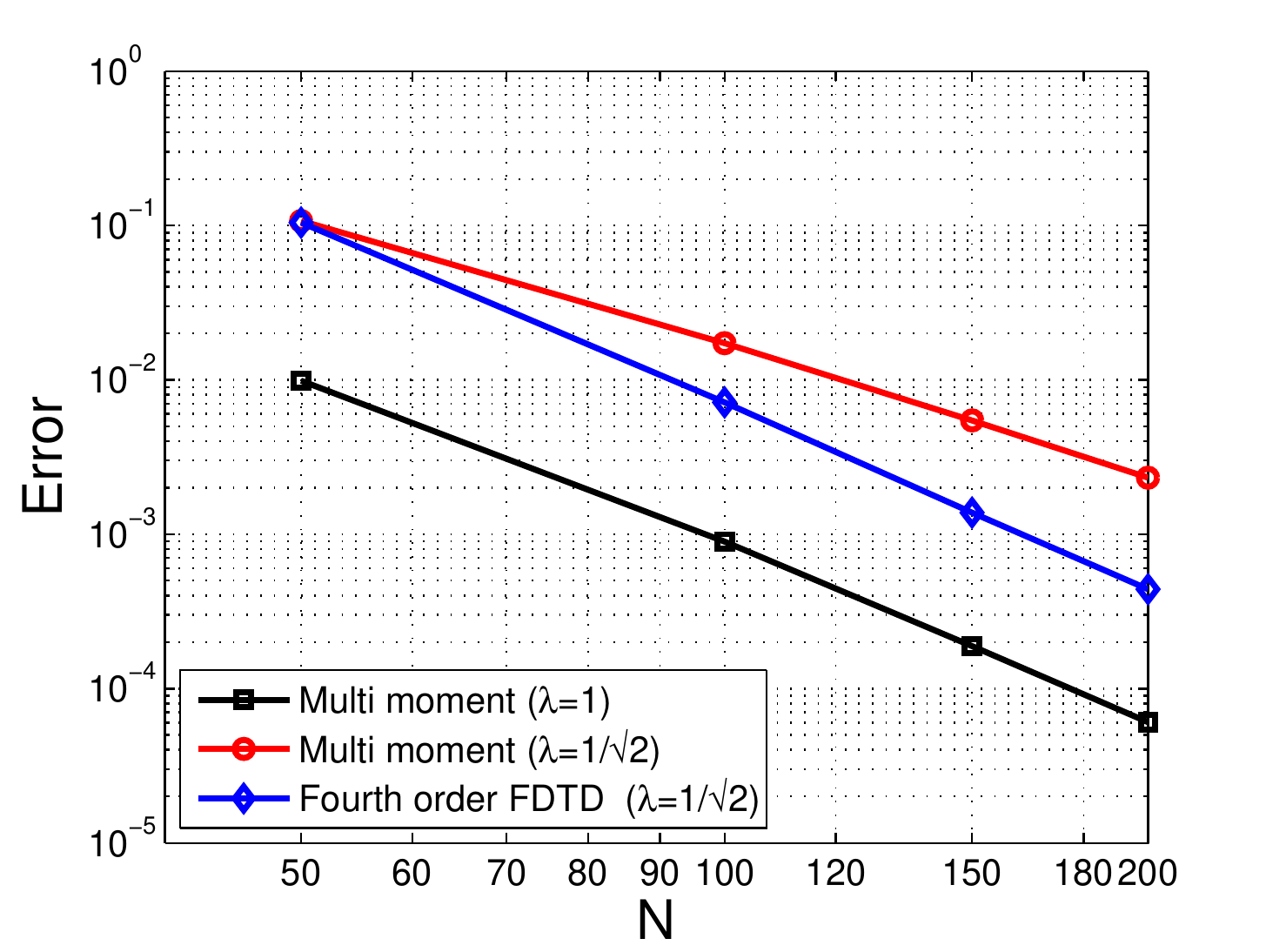}
\includegraphics[width=7.5cm]{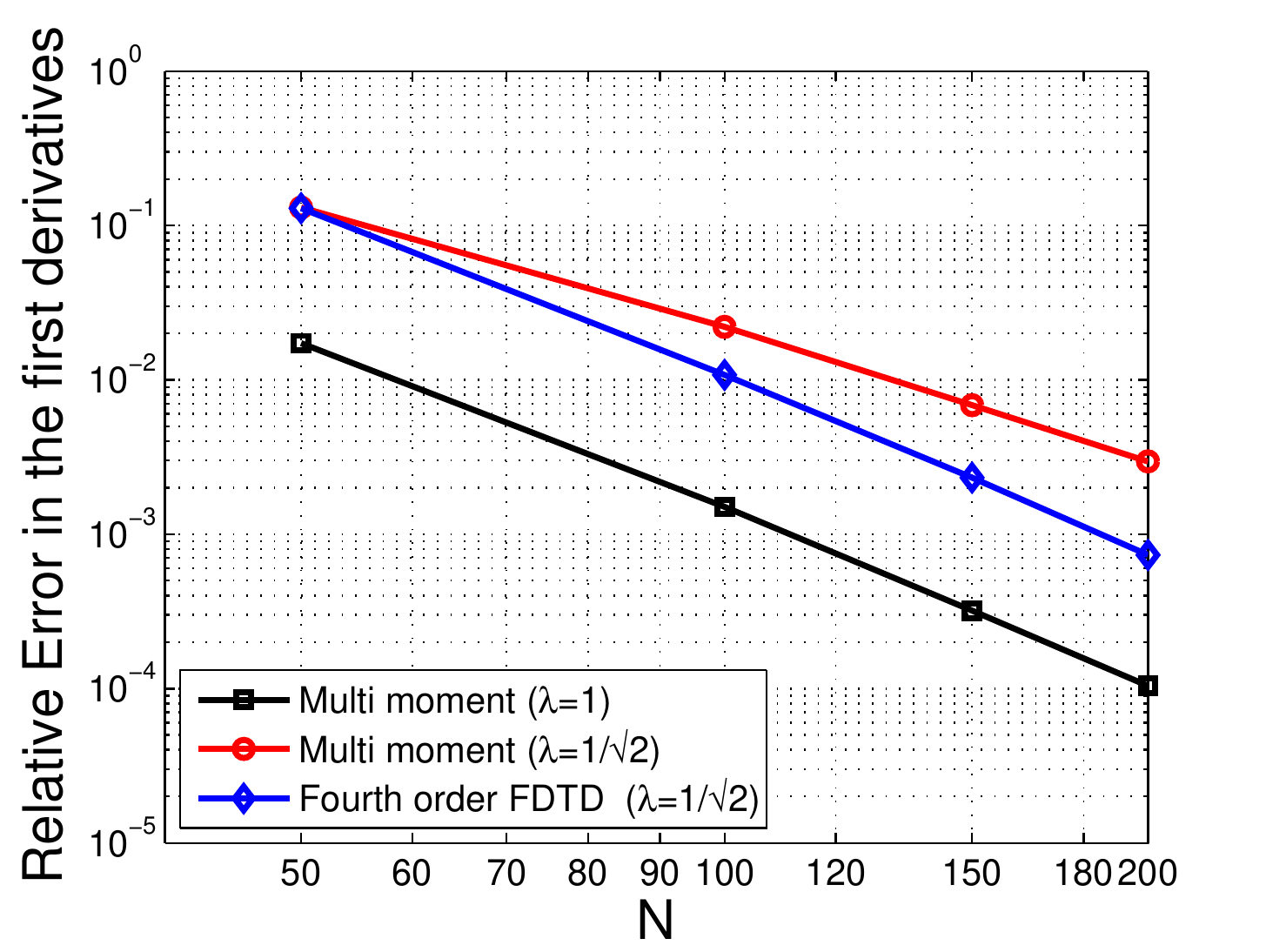}
\caption{Numerical error $\epsilon_1$ (left) and $\epsilon_2$ (right) in the numerical solutions for the initial condition $\sigma^{-2}=500$. The order of accuracy is shown in Table \ref{table:order500}}\label{fig:Error500}
\end{figure}

\begin{figure}[h]
\centering
\includegraphics[width=7.5cm]{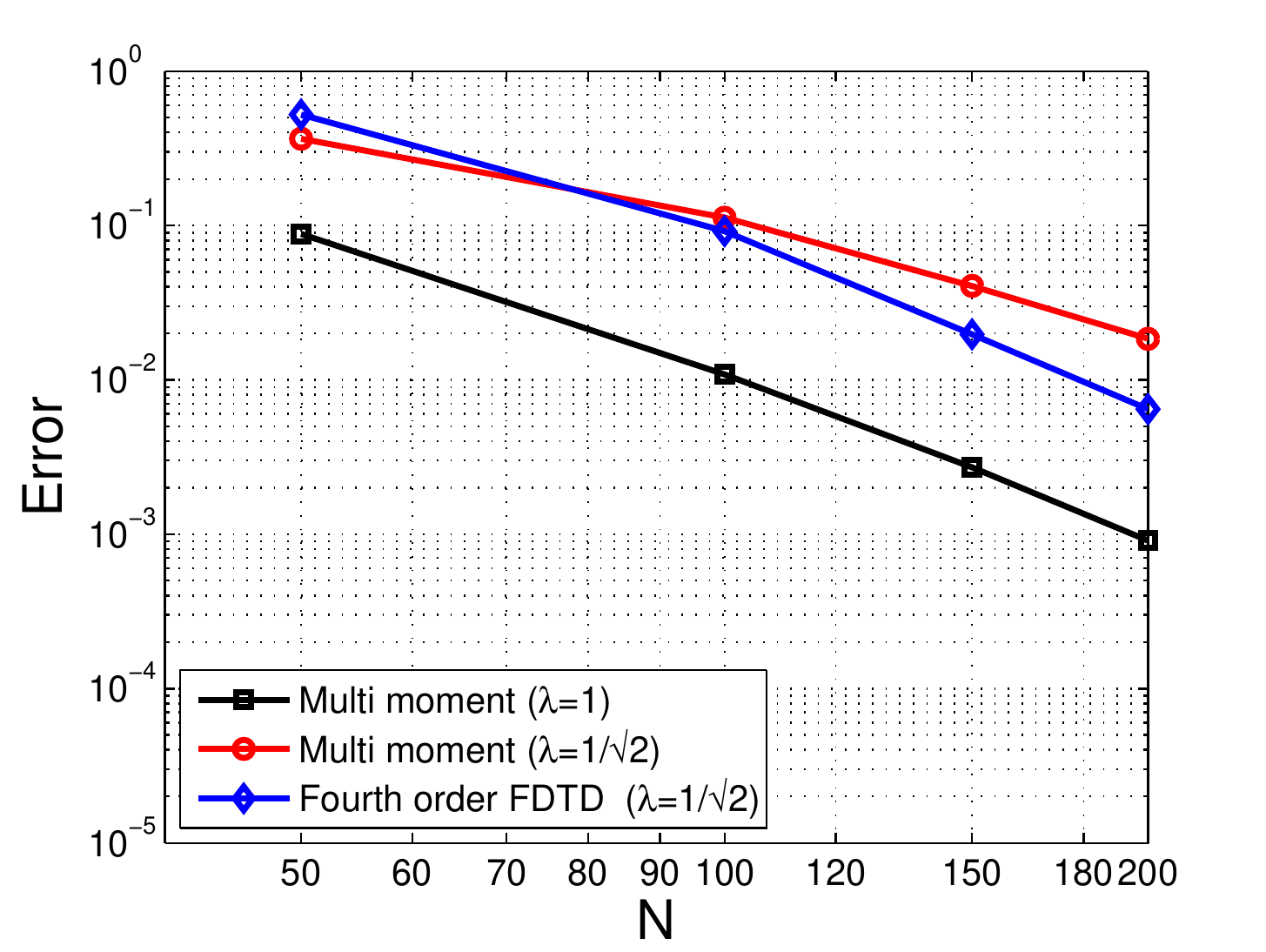}
\includegraphics[width=7.5cm]{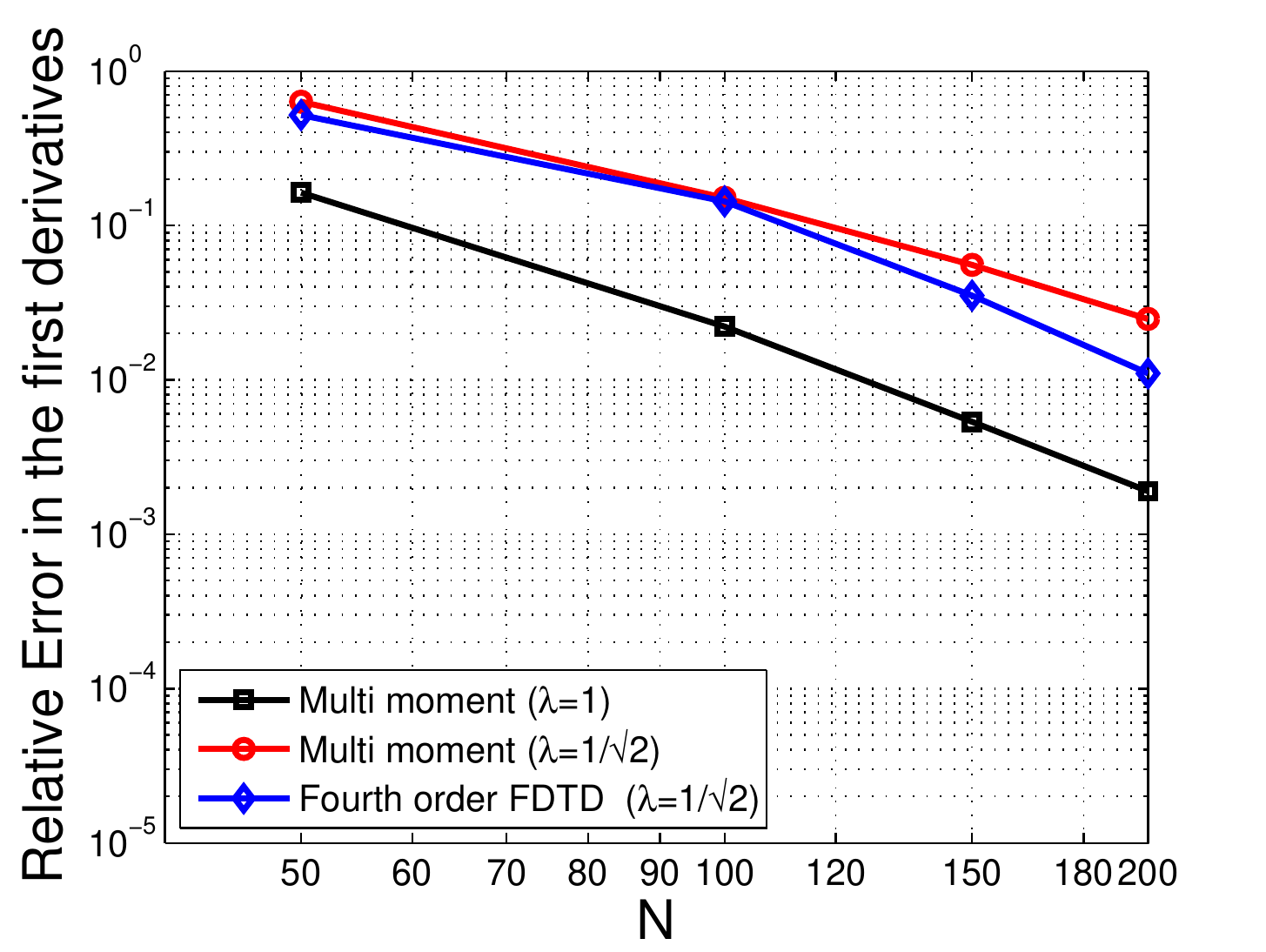}
\caption{Numerical error $\epsilon_1$ (left) and $\epsilon_2$ (right) in the numerical solutions for the initial condition $\sigma^{-2}=1500$. The order of accuracy is shown in Table \ref{table:order1500}}\label{fig:Error1500}
\end{figure}

\begin{table}[ht]
\centering\caption{Numerical errors $\epsilon_1$ and $\epsilon_2$ and the order of convergence. $\sigma$ in the initial profile is $\sigma^{-2}=500$.}
\small\addtolength{\tabcolsep}{-4pt}
\begin{tabular}{c|c|c|c|c|c|c|c|c|c|c|c|c}
  &\multicolumn{4}{|c}{Multi moment ($\lambda=1$)}  &\multicolumn{4}{|c}{Multi moment ($\lambda=\frac{1}{\sqrt{2}}$)}  &\multicolumn{4}{|c}{FDTD ($\lambda=\frac{1}{\sqrt{2}}$)} \\
  \hline
  $N$  & {50}  & 100  & 150  & 200  & {50}  & 100  & 150  & 200  & {50}  & 100  & 150  & 200  \\
  $\epsilon_1$  & 9.86e-3  &8.93e-4  &1.88e-4  &6.06e-5 & 1.07e-1 & 1.72e-2 & 5.45e-3  & 2.32e-3 & 1.04e-1  &7.15e-3 & 1.37e-3 & 4.39e-4\\
  order  & & 3.46 & 3.84 & 3.93 & &2.63 & 2.84 & 2.97  & &  3.86  &4.06 & 3.96 \\
  $\epsilon_2$  & 1.72e-2&  1.50e-3 & 3.20e-4  &1.04e-4 &  1.30e-1&  2.20e-2 & 6.86e-3  &2.95e-3
 & 1.29e-1 & 1.07e-2 & 2.32e-3 & 7.33e-4\\
  order & & 3.51  &3.81  &3.90  & & 2.56 & 2.87 & 2.92 & &3.58 & 3.78  &4.00\\
  \end{tabular} \label{table:order500}
\end{table}

\begin{table}[ht]
\centering\caption{Numerical errors $\epsilon_1$ and $\epsilon_2$ and the order of convergence. $\sigma$ in the initial profile is $\sigma^{-2}=1500$.}\small\addtolength{\tabcolsep}{-4pt}
\begin{tabular}{c|c|c|c|c|c|c|c|c|c|c|c|c}
  &\multicolumn{4}{|c}{Multi moment ($\lambda=1$)}  &\multicolumn{4}{|c}{Multi moment ($\lambda=\frac{1}{\sqrt{2}}$)}  &\multicolumn{4}{|c}{FDTD ($\lambda=\frac{1}{\sqrt{2}}$)} \\
  \hline
  $N$  & {50}  & 100  & 150  & 200  & {50}  & 100  & 150  & 200  & {50}  & 100  & 150  & 200  \\
  $\epsilon_1$ & 8.80e-2  & 1.08e-2  &  2.70e-3 & 9.10e-4  &   3.65e-1 & 1.12e-1  & 4.04e-2 &  1.84e-2  & 5.23e-1  &  9.20e-2 &  1.96e-2&  6.46e-3\\
  order& {}  & 3.02  &  3.41  & 3.79  & {}  &  1.69  & 2.52 & 2.72 &   & 2.50  & 3.80 & 3.86 \\
  $\epsilon_2$ &  1.63e-1  & 2.21e-2  & 5.35e-3  &1.89e-3  &  6.33e-1  & 1.51e-1 & 5.54e-2 & 2.47e-2 &  5.19e-1  &  1.43e-1  & 3.50e-2 & 1.09e-2 \\
  order& {}  &2.87  & 3.50  & 3.60  & {}  & 2.06  & 2.47 & 2.80  &  &  1.85  &3.47& 4.03\\
  \end{tabular} \label{table:order1500}
\end{table}



\begin{figure}[h]
\centering
\includegraphics[scale=0.3]{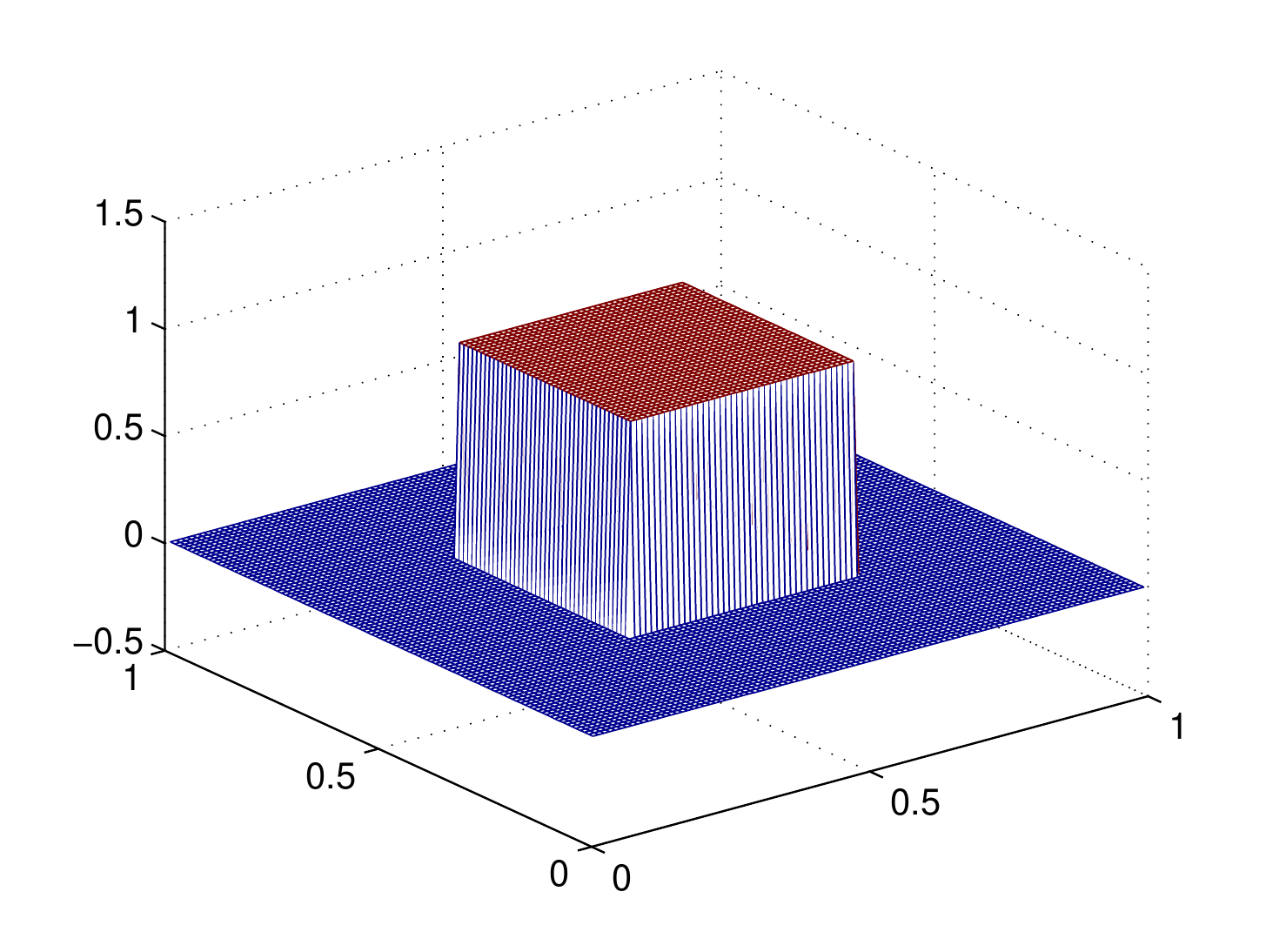}
\includegraphics[scale=0.3]{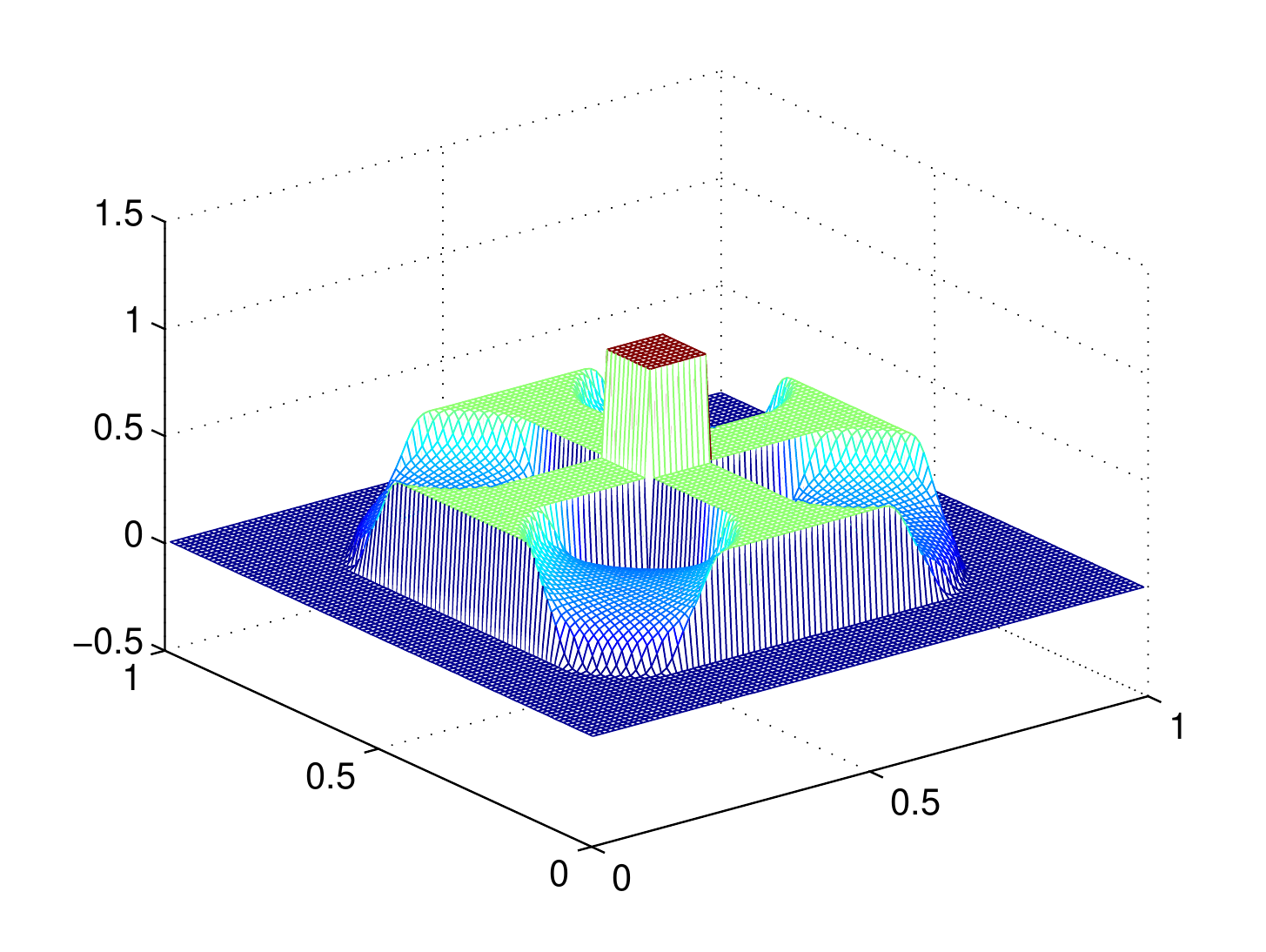}
\includegraphics[scale=0.3]{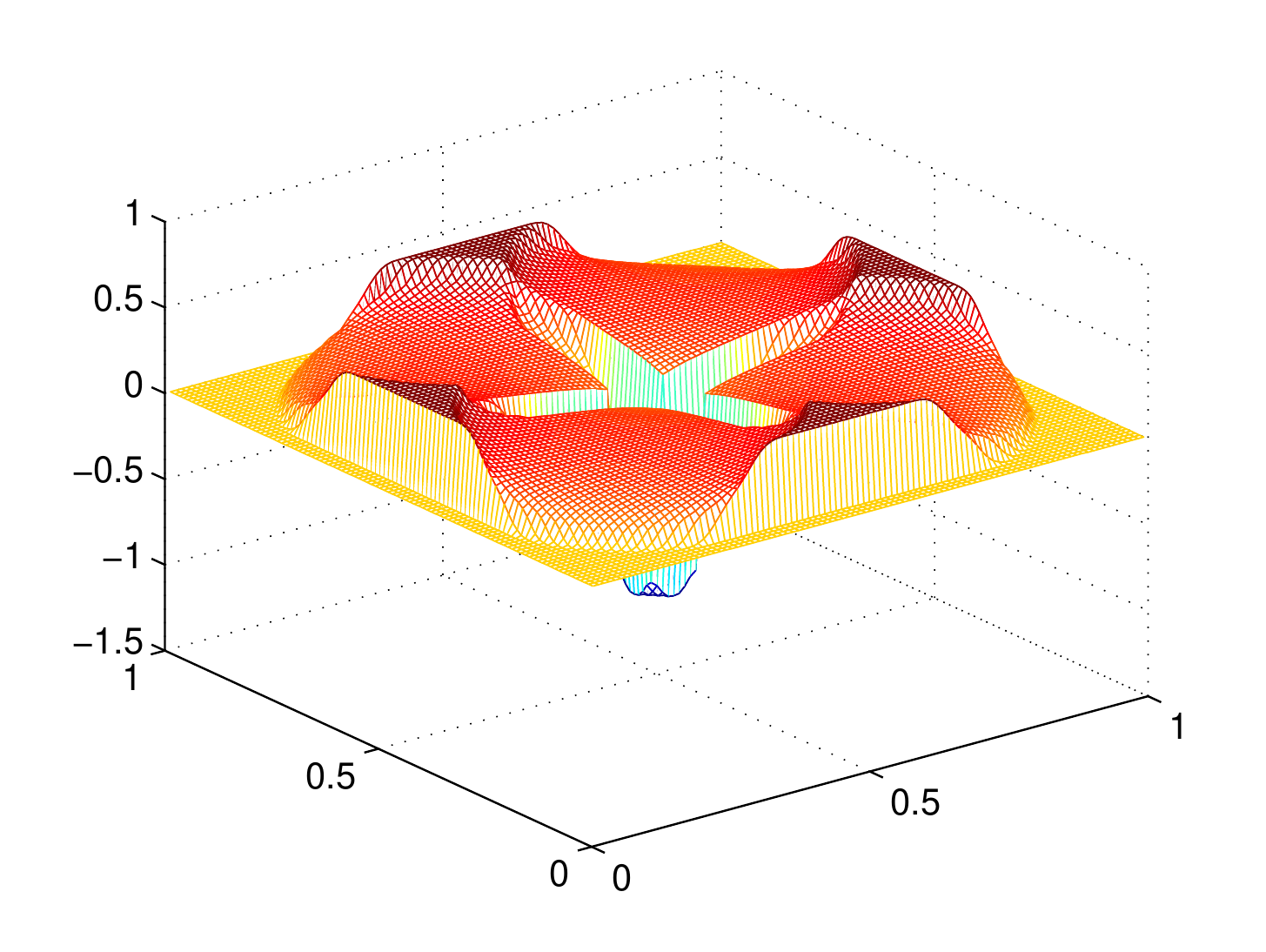}
\includegraphics[scale=0.3]{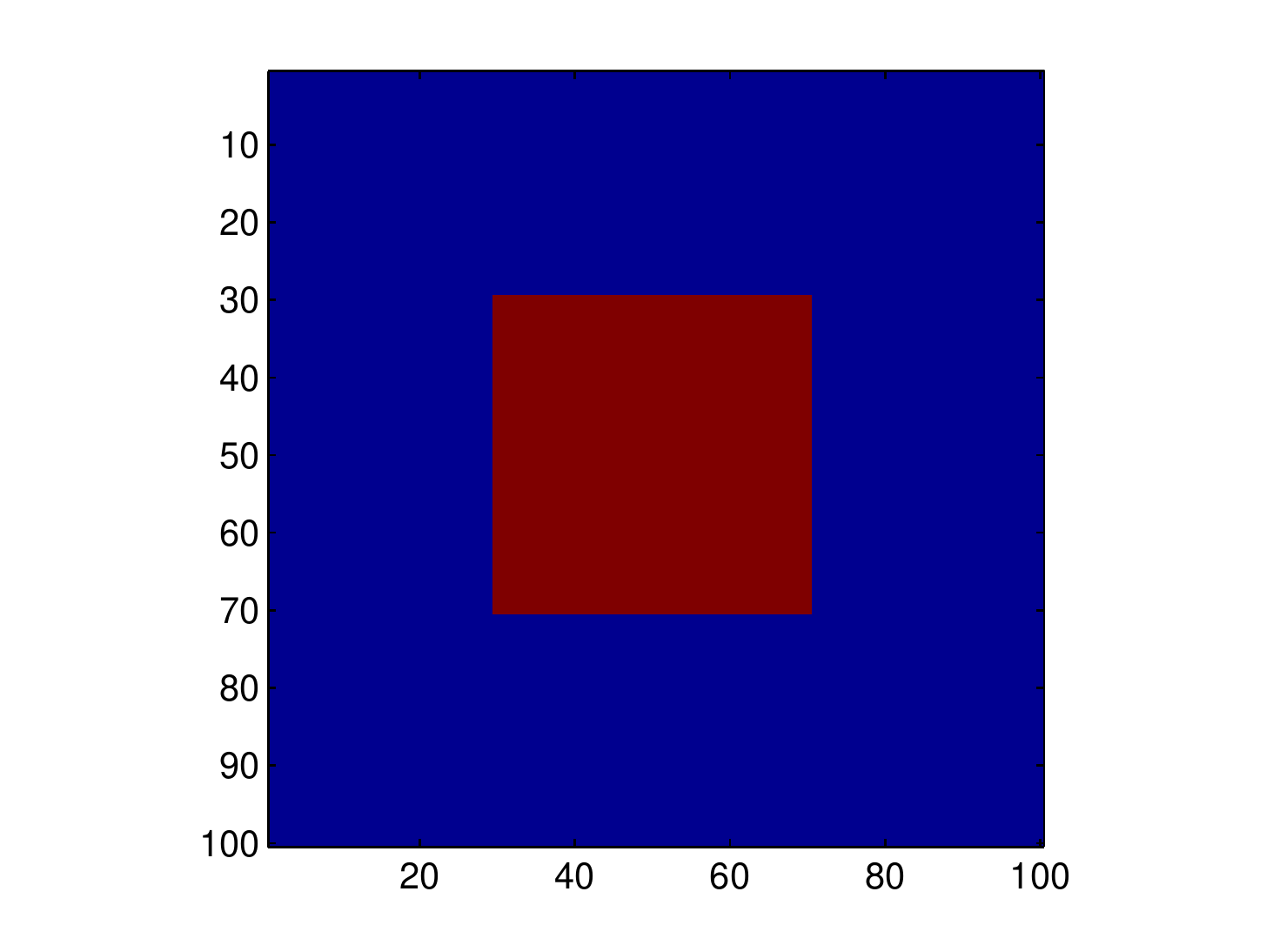}
\includegraphics[scale=0.3]{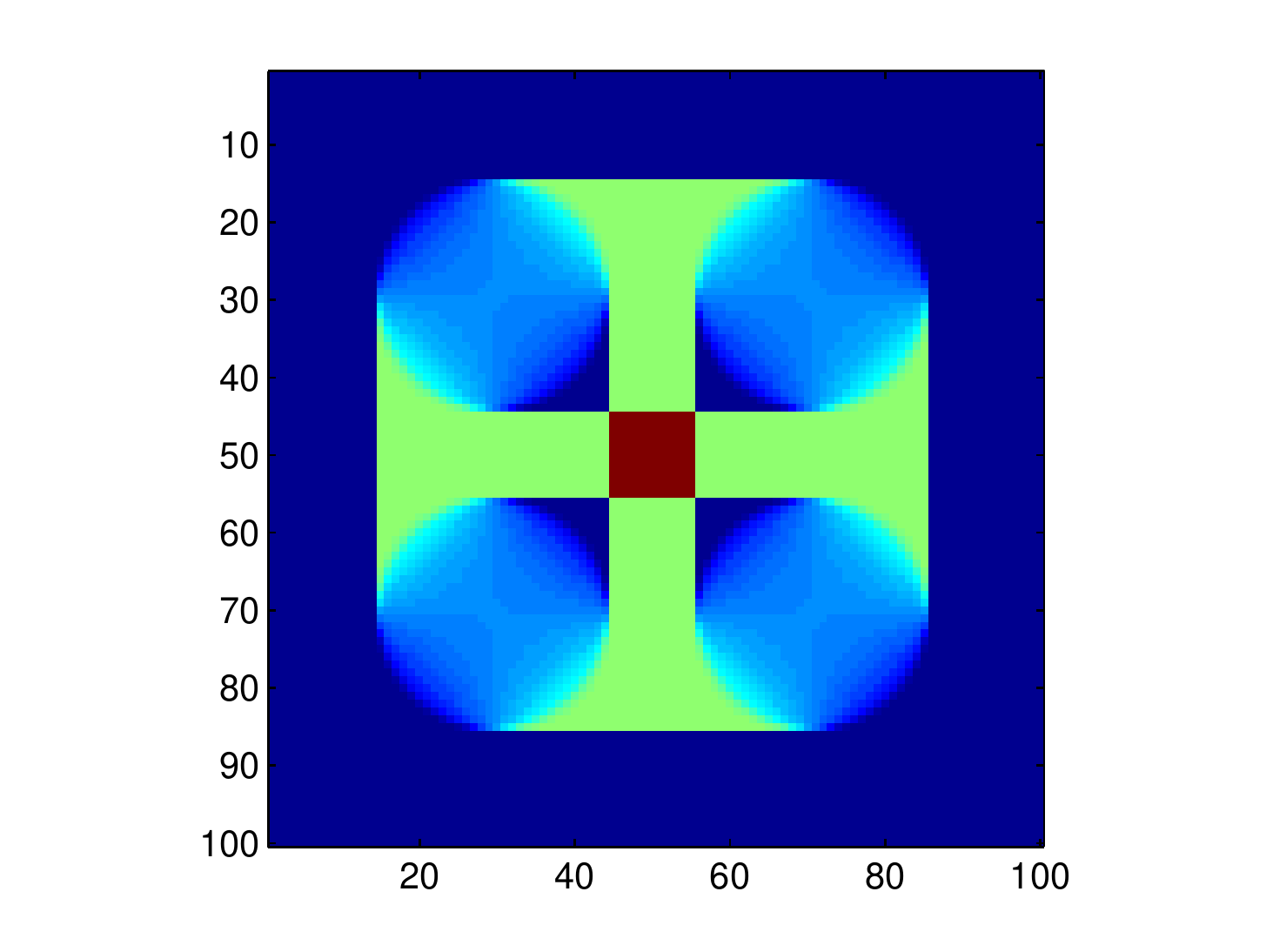}
\includegraphics[scale=0.3]{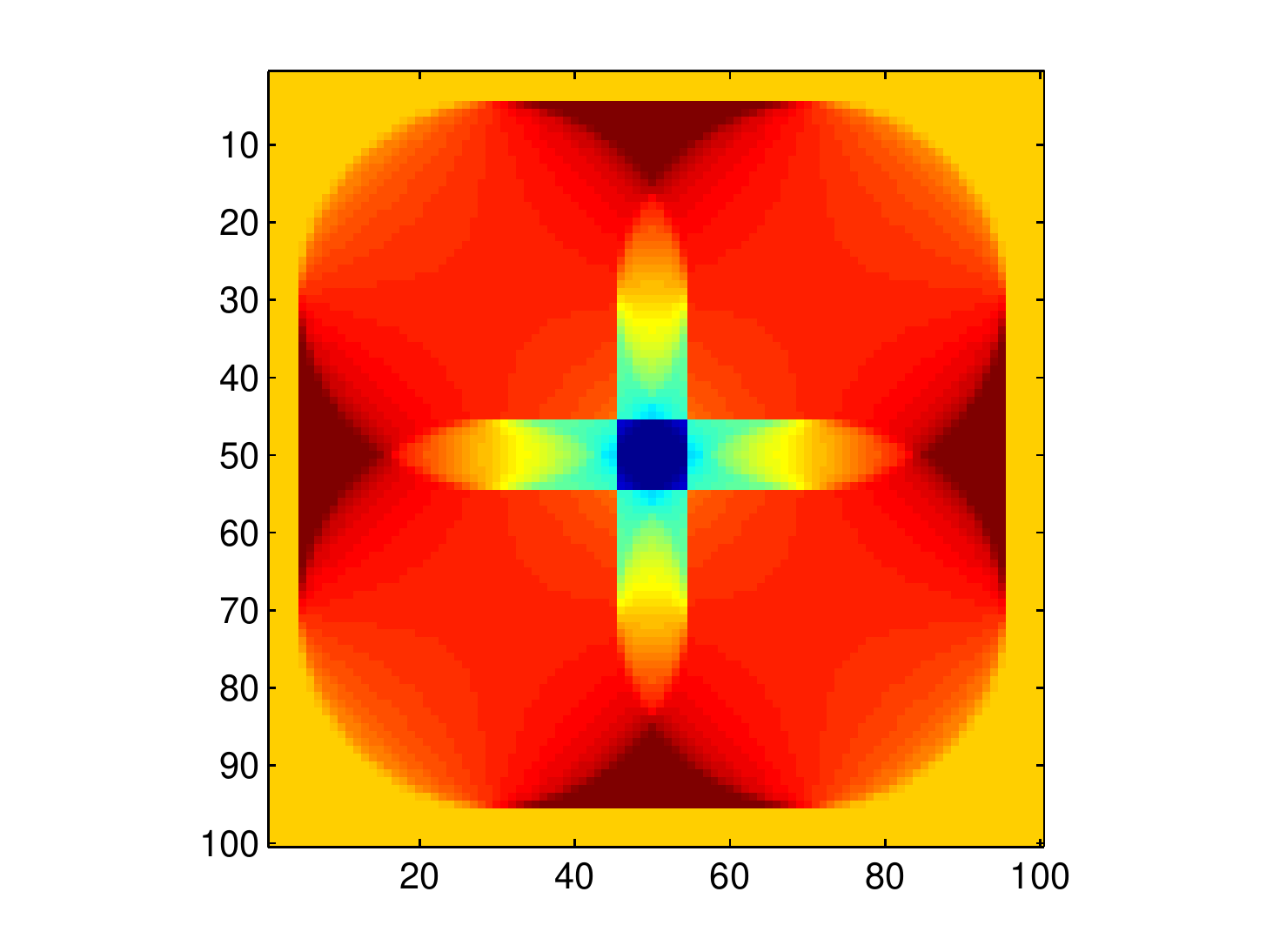}
\caption{The initial condition for $H$ (the first column) and the numerical solution for $H$ at time $T=0.15$ (the second column) and $T=0.25$ (the third column).}\label{fig:square}
\end{figure}

\begin{figure}[h]
\centering
\includegraphics[width=7.8cm]{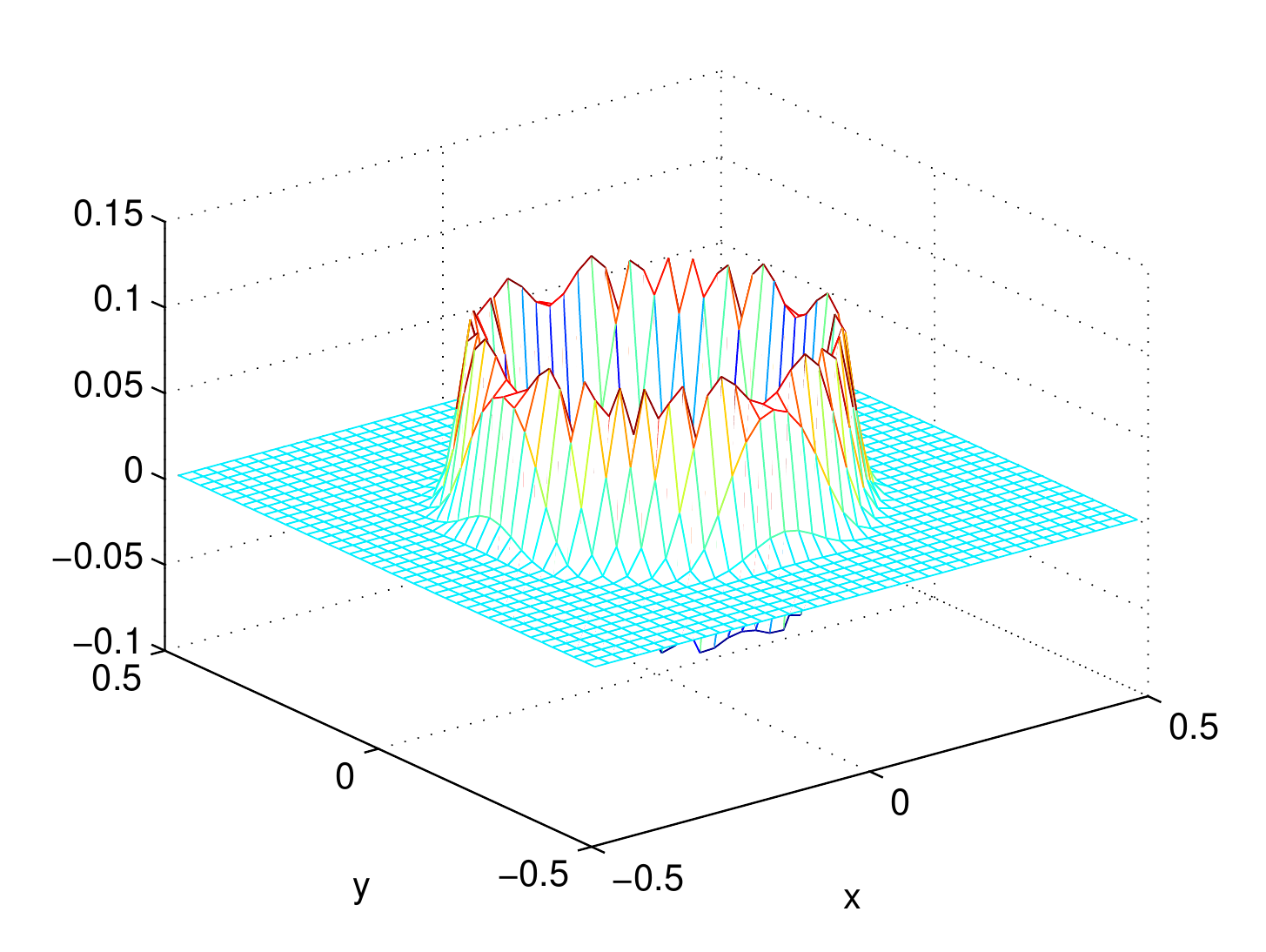}
\includegraphics[width=7.8cm]{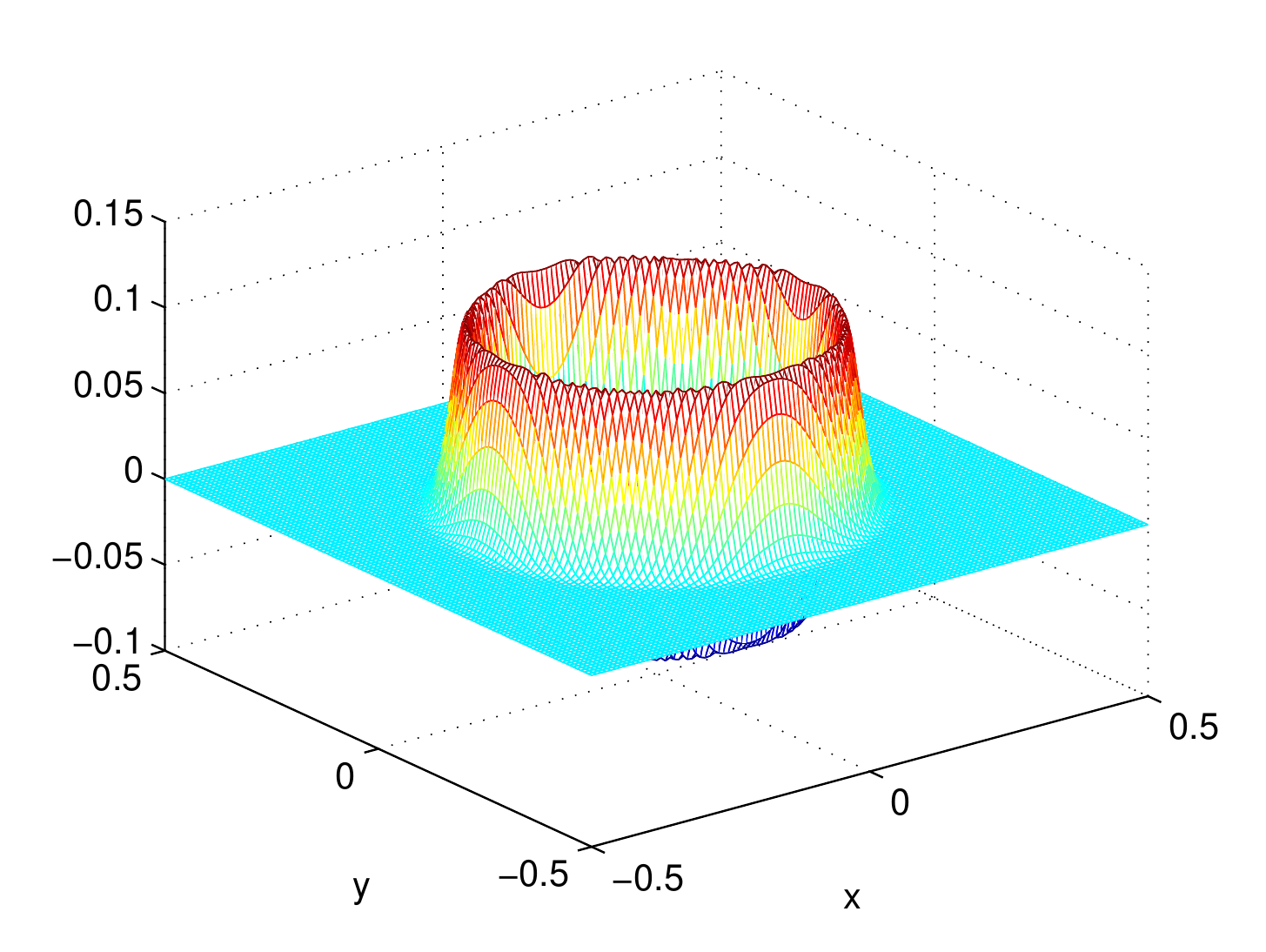}
\includegraphics[width=7.8cm]{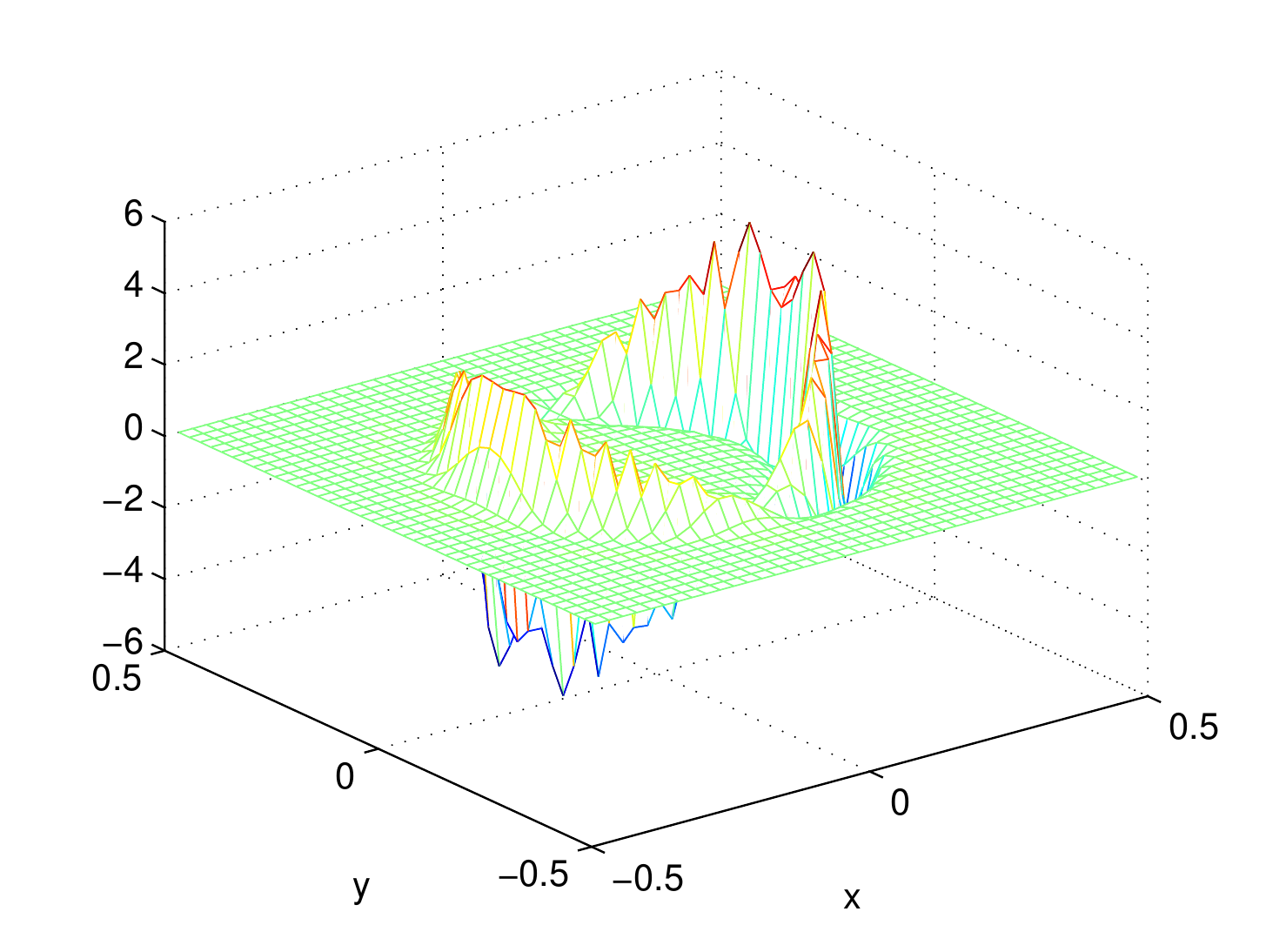}
\includegraphics[width=7.8cm]{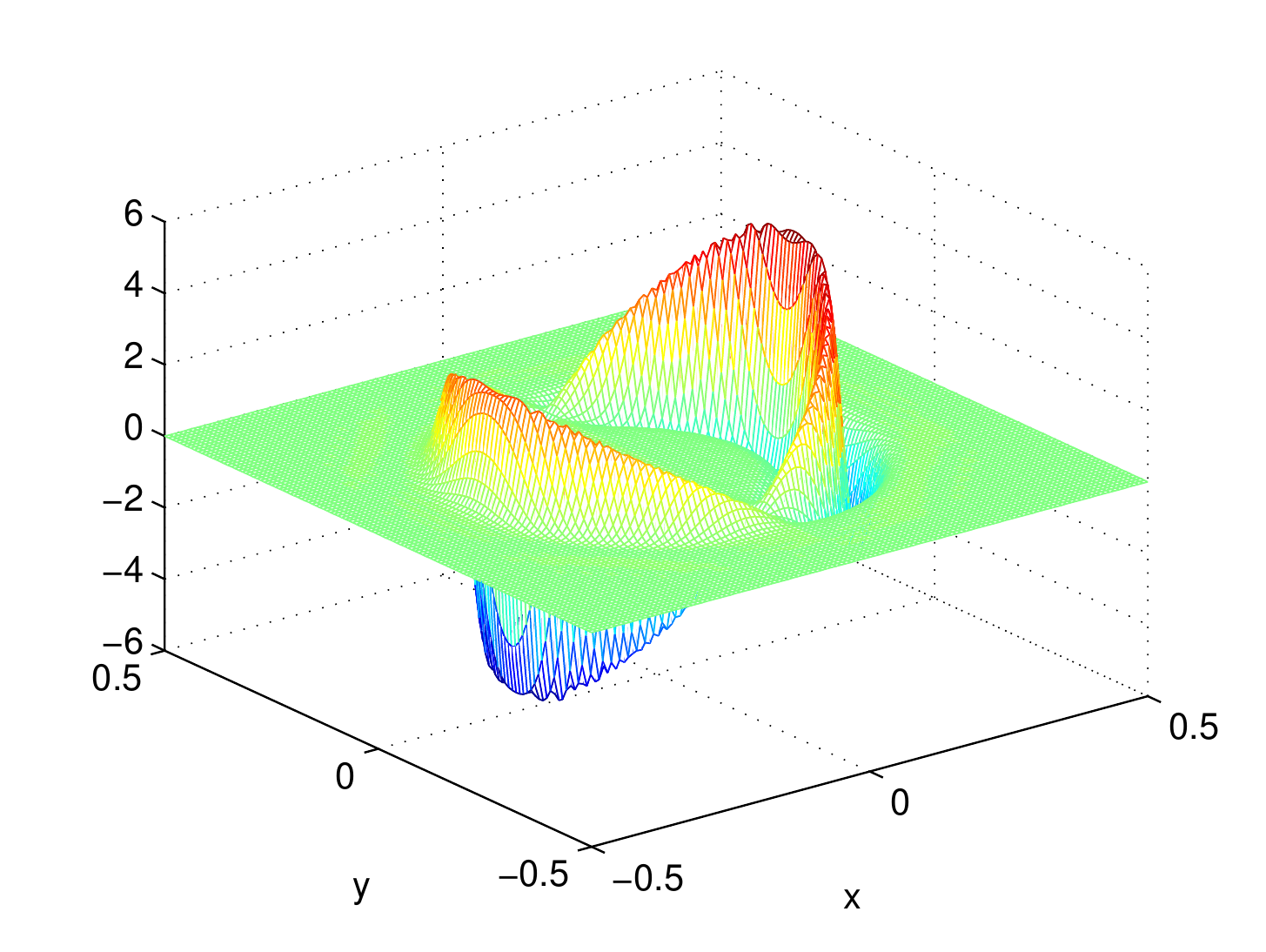}
\caption{Left column: 2D surface plot of $\{ { h^{10}_z}_{i,j} \}_{i,j}$(top) and $\{ {\partial_x h^{10}_z}_{i,j} \}_{i,j}$(bottom). Right column: The bi-cubic interpolation \eqref{bicubicHz} $h_z(x,y)$ (top) and $\partial_xh_z(x,y)$ (bottom). }\label{fig:interp}
\end{figure}
\end{document}